\documentclass{amsart}
\usepackage{mathrsfs}
\usepackage{amsmath}
\usepackage{tikz}
\usepackage{color}
\usepackage{verbatim}
\usepackage[numbers]{natbib}
\usepackage{amscd}
\usepackage{amssymb}
\usepackage{amsthm}
\usepackage{amsfonts}
\usepackage{mathrsfs}
\usepackage{hyperref}
\usepackage{epic}
\usepackage[all]{xy}
\usepackage{mathtools}
\usepackage{enumerate}
\usepackage{geometry}
\newtheorem{Theorem}{Theorem}[section] 
\newtheorem{Lemma}[Theorem]{Lemma}
\newtheorem{Proposition}[Theorem]{Proposition}
\newtheorem{Corollary}[Theorem]{Corollary}

\begin{document}
\setlength{\unitlength}{1mm}
\title{Non-trivial Kazhdan-Lusztig coefficients of affine Weyl groups}
\author{Pan Chen}
\address{Faculty of Mathematics\\            
Yunnan Normal University\\
Chenggong District\\
KunMing, 650500\\
P. R. China}
\email{chenpan901@163.com }
\maketitle

%

\begin{abstract}
In this paper we show that the leading coefficient $\mu(y,w)$ of some Kazhdan-Lusztig polynomials $P_{y,w}$
with $y,w$ in an affine Weyl group of type $\tilde B_n$ (resp. $\tilde C_n$ or $\tilde D_n$)  is  $n$ (resp. $n+1$).
\end{abstract}

\section{Introduction}

In Kazhdan-Lusztig theory, an interesting question is to compute the Kazhdan-Lusztig coefficient $\mu(y,w)$ of an Kazhdan-Lusztig polynomial $P_{y, w}$.
Even for symmetric groups, the coefficient $\mu(y,w)$ is far from being understood. McLarnan and Warrington \cite{mclarnan2003counter} showed that for symmetric group $\mathfrak S_{10}$, $\mu(y,w)$ can be greater than 1, so the (0,1)-conjecture has a negative answer. Further, Warrington \cite{warrington2011equiv} shows that the upper bound of  $\mu(y,w)$ for symmetric group  $\mathfrak S_{n}$ should increase rapidly when  $n$ increases.

For affine Weyl groups, also not much work is done.
Lusztig \cite{lusztig1997nonlocal} computed the  coefficient $\mu(y,w)$  of some Kazhdan-Lusztig polynomials of an affine Weyl group of type
$\tilde B_2$ ; more were computed in \cite{wang2008leading} and \cite{fwang2021leading}. Scott \cite{scott2003some} found some non-trivial examples  of the coefficient (larger than 1)
for an affine Weyl group of type $\tilde A_5$. Xi \cite{xi2005leading} showed that $\mu(y,w)\leqslant 1$ if $a(y)<a(w)$
when $W$ is an affine Weyl group of type $\tilde A_n$; for the same affine Weyl group, Green \cite{green2009leading}  proved that $\mu(y,w)\leqslant 1$
if $y$ is a fully commutative element.
In \cite{scott2010some} Scott and Xi showed that some $\mu(y,w)$ is $n+2$ for an
affine Weyl group of type $\tilde A_n$.
Guo and Zheng \cite{gz2021leading} showed that $\mu(y,w)\leq 3$ if $y,w$ are in the lowest two-sided cell of the affine Weyl group of type $\tilde G_2$.

Motivated by \cite{scott2010some}, in this paper we show that the  coefficient $\mu(y,w)$ of some Kazhdan-Lusztig polynomials $P_{y,w}$
with $y,w$ in an affine Weyl group of type $\tilde B_n$ (resp. $\tilde C_n$ or $\tilde D_n$)  is  $n$ (resp. $n+1$), see Theorem \ref{thm:Bmu}, Theorem \ref{thm:Cmu}
and Theorem \ref{thm:Dmu}.

From
\cite{andersen1986inversion} and \cite{andersen1994representations} we know that some of the  coefficients $\mu(y,w)$ for an affine Weyl group can be identified with the dimension of
first extension group of some
irreducible
modules of the corresponding algebraic group over $\bar{\mathbb{F}}_p$, provided that $p$ is sufficiently large. This fact  also increases our interest in computing the coefficient for affine Weyl groups.

This paper is organized as follows. In section 2, we recall some basic facts about a formula of Springer and about the lowest
two-sided cell of an affine Weyl group. In  sections 3, 4 and 5, we deal with the cases of types $\tilde B_n,\ \tilde C_n$
and $\tilde D_n$  respectively. Our main results are Theorem \ref{thm:Bmu}, Theorem \ref{thm:Cmu}
and Theorem \ref{thm:Dmu} respectively, which say that the coefficient $\mu(y,w)$ is $n,\ n+1,\ n+1$ respectively for some $y,\ w$ in an affine Weyl group of  type $\tilde B_n,\ \tilde C_n,\ \tilde D_n$ respectively.  The key points are   to find the right elements $y,w$ and then to compute   value of $\mu(y,w)$.

\section{Preliminaries}

In this section we recall some facts and fix some notations. Let $G$ be a connected, simply connected simple algebraic group over the field $\mathbb{C}$ of complex numbers and let $T$ be a maximal torus
of $G$. Then the Weyl group $W_0=N_G(T)/T$ acts naturally on the character group $X=\text{Hom}(T,\mathbb{C}^*)$ of $T$. The semidirect product
$W=W_0\ltimes X$ is called an extended affine Weyl group. It contains the affine Weyl group $W_a=W_0\ltimes \mathbb{Z}R$, where $R$ is the root
system of
$G$ with respect to $T$. For $x\in X,$ denote by $t_x$ the translation $y\mapsto y+x$ of $X.$ Let $E=\mathbb{R}\otimes X.$ Then $W$ acts on $E$ naturally and the action is faithful. Thus $W$ can be identified with a subgroup of the affine transformation
group on $E.$

Fix a simple root system $\Delta\subset R$. Let $\alpha_0$ be the highest short root of $R$. Let $s_\alpha$ be the reflection corresponding to the root $\alpha\in R.$ Set
$s_0=s_{\alpha_0}t_{\alpha_0},S=\{s_\alpha,s_0|\alpha\in \Delta\},$ then $(W_a,S)$ is a Coxeter group.

We have $W=W_a\rtimes \Omega,$ where $\Omega=\{u\in W|uSu^{-1}=S\}.$ For any $\alpha\in R,\ k\in \mathbb{Z},$ set $s_{\alpha,k}(y)=y-(\langle
y,\alpha^{\vee}\rangle-k)\alpha,\ y\in E.$ Then we have $wt_z w^{-1}=t_{w(z)}$, $s_{\alpha,k}=t_{k\alpha}s_\alpha=s_\alpha t_{-k\alpha}$ and
$(wt_z)s_{\alpha,k}(wt_z)^{-1}=s_{w(\alpha),k+\langle z,\alpha^{\vee}\rangle}$ for any $z\in X,\ \alpha\in R,\ w\in W_0,\ k\in \mathbb{Z}.$

Let $l$ and $\leqslant$ be the length function and the Bruhat order of $(W_a,S)$ respectively.  Lusztig\cite{Lusztig1983Singularities} has extended them to $W$ as follows: $l(\omega y)=l(y)$, and  $\omega y\le \omega' w$ if and only if $y\le w$ and $\omega=\omega'$, where $y,\ w\in W_a,\ \omega,\omega'\in \Omega$. Then the Hecke algebra $H$ over $\mathcal{A}=\mathbb{Z}[q^{\frac{1}{2}},q^{-\frac{1}{2}}]$ of  $(W,S)$ is defined to be the free $\mathcal A$-module with a basis $T_w, \ w\in W$ and multiplication relations
\begin{alignat*}{2} (T_s-q)(T_s+1)&=0,\qquad s\in S;\\
T_wT_u&=T_{wu},\qquad w,u\in W, \ l(wu)=l(w)+l(u).\end{alignat*}
So Kazhdan-Lusztig polynomial $P_{y,w}$ for any $y,\ w\in W$ is defined. We shall use the Kazhdan-Lusztig basis
$C_w=q^{-\frac{l(w)}{2}}\sum\limits_{y\leqslant w}P_{y,w}T_y, w\in W$, of the Hecke algebra $H$ of $(W,S)$.

Recall that  $P_{y,w}=\mu(y,w)q^{\frac{1}{2}(l(w)-l(y)-1)}+\text{lower degree terms}$ if $y<w$. This defines  the Kazhdan-Lusztig coefficient $\mu(y,w)$ of $P_{y,w}$. If $w<y$, set $\mu(y,w)=\mu(w,y).$  If $y\not < w$ and $w\not < y$, set $\mu(y,w)=0.$

If $y<w$ and $\mu(y,w)\neq 0$  we write $y\prec w.$ Denote $\mathcal{L}(w)=\{s\in S|sw<w\}.$ We collect some well-known properties about $P_{y,w}$ and $\mu(y,w)$.

\begin{Proposition}\label{prop:kl}(\cite{kazhdan1979representations}) Assume that $y,w\in W_a$ and $y\leqslant w.$
\begin{enumerate}[(1)]
    \item If $sw<w$ for some $s\in S,$  then
$$P_{y,w}=q^{1-c}P_{sy,sw}+q^cP_{y,sw}-\sum\limits_{\mbox{\tiny$\begin{array}{c} y\leqslant z\prec sw\\ sz<z\end{array}$}}\mu(z,sw)q^{\frac{1}{2}(l(w)-l(z))}P_{y,z},$$
where $c=1$ if $sy<y,$ and $c=0$ if $sy>y.$ (Convention: $P_{y,w}=0$ if $y\not\leqslant w.$)\\
    \item If  $s\in S$ and $sw<w$ (resp. $ws<w$), then $P_{y,w}=P_{sy,w}$ (resp. $P_{y,w}=P_{ys,w}$).\\
    \item If  $s\in S$ and $sw<w,\ y\not\leqslant sw,$ then
$$P_{y,w}=P_{sy,sw}.$$
    \item If $s\in S$ and $sy>y$, $sw<w$, then $y\prec w$ if and only if $w=sy.$ Furthermore, in this case $\mu(y,w)=1.$\\
    \item Let $s,t\in S$ be such that  $st$ has order 3. Set $\mathscr{D}_{\mathcal{L}}(s,t)=\{w\in W_a|\#\mathcal{L}(w)\cap\{s,t\}=1\}.$ If $v\in \mathscr{D}_{\mathcal{L}}(s,t),$ then exactly one of the elements $sv,tv$ is in $\mathscr{D}_{\mathcal{L}}(s,t),$ denoted by ${}^*v.$ The map $v\mapsto {}^*v$ is an involution of $\mathscr{D}_{\mathcal{L}}(s,t),$ called the (left) star operator. For $y,w\in \mathscr{D}_{\mathcal{L}}(s,t),$ we have
$$\mu(y,w)=\mu({}^*y, {}^*w).$$
\end{enumerate}
\end{Proposition}

From the above Proposition \ref{prop:kl} (4) we know that if $y\prec w$ and $l(w)-l(y)>1$ then
$\mathcal{L}(y)\supseteq \mathcal{L}(w)$ and $l(w)-l(y)\equiv 1 (\text{mod}\ 2)$. This fact will be used frequently in our calculations.

Write
$$C_xC_y=\sum\limits_{z\in W}h_{x,y,z}C_z,\quad  h_{x,y,z}\in \mathcal{A}.$$
Denote the highest degree of $h_{x,y,z}$ in the invariable $q^{\frac{1}{2}}$ by $\deg_{q^{\frac{1}{2}}}h_{x,y,z}$. The degree is controlled by the length
function.

\begin{Proposition}\label{prop:length}
$\deg_{q^{\frac{1}{2}}}h_{x,y,z}\leqslant l(x)+l(y)-l(z).$
\end{Proposition}

\begin{proof}
We use induction on $l(x)\in \mathbb{N}.$ When $l(x)=0$, it's obvious. When $l(x)>0$, there exists $s\in S,$ such that $xs<x.$ Let
$x^\prime=xs.$
Now
$$C_xC_y=(C_{x^\prime}C_s-\sum\limits_{vs<v\prec x'}\mu(v,x')C_v)C_y=\sum\limits_{z\in W}h_{x,y,z}C_z.$$
There are two cases:
\\(a) $sy<y.$

In this case, $C_sC_y=(q^{\frac{1}{2}}+q^{-\frac{1}{2}})C_y.$ Thus
$$h_{x,y,z}=(q^{\frac{1}{2}}+q^{-\frac{1}{2}})h_{x^\prime,y,z}-\sum\limits_{vs<v\prec x'}\mu(v,x')h_{v,y,z},$$
and
$$\deg_{q^{\frac{1}{2}}}h_{x,y,z}\leqslant \max\limits_{
vs<v\prec x'}\{\deg_{q^{\frac{1}{2}}}h_{x^\prime,y,z}+1,\deg_{q^{\frac{1}{2}}}h_{v,y,z}\}.$$
The conclusion follows from the induction hypothesis.
\\(b) $sy>y.$

In this case, $C_sC_y=C_{sy}+\sum\limits_{sw<w\prec y}\mu(w,y)C_w.$ Thus
$$h_{x,y,z}=h_{x^\prime,sy,z}+\sum\limits_{sw<w\prec y}\mu(w,y)h_{x^\prime,w,z}-\sum\limits_{vs<v\prec x'}\mu(v,x')h_{v,y,z},$$
and
$$\deg_{q^{\frac{1}{2}}}h_{x,y,z}\leqslant \max\limits_{\mbox{\tiny
$\begin{array}{c} sw<w\prec y\\ vs<v\prec x'
\end{array}$}}\{\deg_{q^{\frac{1}{2}}}h_{x^\prime,sy,z},\deg_{q^{\frac{1}{2}}}h_{x^\prime,w,z},\deg_{q^{\frac{1}{2}}}h_{v,y,z}\}.$$
By the induction hypothesis, we have the needed inequality.
\end{proof}

Following Lusztig, define
$$a(z)=\min\{i\in \mathbb{N}|q^{-\frac{i}{2}}h_{x,y,z}\in \mathbb{Z}[q^{-\frac{1}{2}}]\text{ for all }x,y\in W\}.$$
If for any $i$ there exist $x,y\in W$ such that $q^{-\frac{i}{2}}h_{x,y,z}\notin \mathbb{Z}[q^{-\frac{1}{2}}]$, we set $a(z)=\infty$. Then $a(w)\leqslant
l(w_0)$ for all $w\in W$, where $w_0$ is the longest element of $W_0$ (see \citep{lusztig1985cells}). Write
$$h_{x,y,z}=\gamma_{x,y,z}q^{\frac{a(z)}{2}}+\delta_{x,y,z}q^{\frac{a(z)-1}{2}}+\text{lower degree terms}.$$
Let $e$ be the neutral element of $W$. Springer showed that $l(z)-a(z)-2\delta(z)\geqslant 0$, where $\delta(z)$ is the degree of $P_{e,z}$
(see \citep{lusztig1987cells}).

Set
$$\mathcal{D}_i=\{z\in W|l(z)-a(z)-2\delta(z)=i\}.$$
Let $P_{e,z}=\pi(z)q^{\delta(z)}$+lower degree terms. The following formula is due to Springer, and the proof can been found in \citep{xi2005leading}.

\begin{Proposition}\label{prop:Springer}
(Springer's formula) For any $x,y\in W,$
$$\mu(y,x)=\sum\limits_{d\in \mathcal{D}_0}\delta_{y^{-1},x,d}+\sum\limits_{f\in \mathcal{D}_1}\gamma_{y^{-1},x,f}\pi(f).$$
\end{Proposition}

We refer to \citep{kazhdan1979representations} and \citep{Lusztig1983Singularities} for the definitions of the preorders $\leqslant _L, \leqslant_R, \leqslant_{LR}$ and of the
equivalence relations $\sim_L,\sim_R,\sim_{LR}$ on $W$. The corresponding equivalence classes are called left cells, right cells, two-sided
cells of $W$ respectively. The preorder $\leqslant_L$ (resp. $\leqslant_R; \leqslant_{LR}$) induces a partial order on the set of left (resp.
right; two-sided) cells of $W$. Denote the induced partial order again by $\leqslant _L$ (resp. $ \leqslant_R; \leqslant_{LR}$). We list some well-known properties.

\begin{Proposition}\label{prop:cell}
(see \citep{lusztig1985cells},\citep{lusztig1987cells} and \citep{xi2005leading})
\begin{enumerate}[(1)]
    \item If $h_{x,y,z}\neq 0,$ then $z\leqslant_R x,\ z\leqslant_L y.$
    \item If $\gamma_{x,y,z}\neq 0,$ then $x\sim_L y^{-1},\ y\sim_L z,\ z\sim_R x.$
    \item If $\delta_{x,y,z}\neq 0,$ then $z\sim_L y$ or $z\sim_R x. $
    \item If $x\sim_{LR} y,$ and $x\leqslant_L y$ (resp. $x\leqslant_R y$), then $x\sim_L y$ (resp. $x\sim_R y$).
\end{enumerate}
\end{Proposition}

Recall that $w_0$ is the longest element of $W_0.$ It is known that $c_0=\{w\in W|a(w)=l(w_0)\}$ is a two-sided cell. In fact, it is the lowest
two-sided cell with respect to the partial order $\leqslant_{LR}.$ See \citep{shi1987two} and \citep{lusztig1985cells}.

Let $R^+$ (resp. $R^-;\ \Delta$) be the set of positive (resp. negative; simple) roots of the root system $R$.  Set $X^+=\{x\in X|\langle
x,\alpha^{\vee}\rangle\geqslant 0 \text{ for all }\alpha\in \Delta\}$ to be  the set of dominant weights of $X.$ Then we have the length formula
(see \citep{iwahori1965some})
$$l(xw)=\sum\limits_{\mbox{\tiny$\begin{array}{c} \alpha\in R^+\\ w^{-1}(\alpha)\in R^- \end{array}$}}|\langle
x,\alpha^{\vee}\rangle+1|+\sum\limits_{\mbox{\tiny$\begin{array}{c} \alpha\in R^+\\ w^{-1}(\alpha)\in R^+ \end{array}$}}|\langle
x,\alpha^{\vee}\rangle|,\text{ for }w\in W_0, x\in X.$$
Thus
$$l(x)=\langle x,2\rho^{\vee}\rangle,\text{ for any }x\in X^+,\text{ where }\rho^{\vee}=\frac{1}{2}\sum\limits_{\alpha\in R^+}\alpha^{\vee}.$$

For each simple root $\alpha$ there is a corresponding fundamental weight $x_\alpha$ in $X$ such that $\langle
x_\alpha,\beta^{\vee}\rangle=\delta_{\alpha,\beta}$ for any simple root $\beta$. For each $w\in W_0,$ set
$$d_w=w\prod\limits_{\mbox{\tiny $\begin{array}{c} \alpha\in\Delta\\ w(\alpha)\in R^- \end{array}$}}x_\alpha.$$
Then
$$c_0=\{d_wxw_0d_u^{-1}|w,u\in W_0,x\in X^+\}.$$
Moreover, for any $w,u\in W_0$, the set $c_{0,w}^\prime=\{d_wxw_0d_u^{-1}|u\in W_0,x\in X^+\}$ is a right cell of $W$ and the set
$c_{0,u}=\{d_wxw_0d_u^{-1}|w\in W_0,x\in X^+\}$ is a left cell of $W$ (see \citep{jian1988two} and \citep{scott2010some}). We shall denote the left cell $c_{0,e}$ by
$\Gamma_0.$  We have (see \citep{lusztig1985cells})
$$\Gamma_0=\{ww_0|w\in W,l(ww_0)=l(w)+l(w_0)\}.$$
The following is a useful property about $\Gamma_0.$

\begin{Proposition} \label{prop:leftcell}
Let $x,y\in W_a$ be such that $l(xw_0)=l(x)+l(w_0),\ l(yw_0)=l(y)+l(w_0).$ Then $xw_0\leqslant yw_0$ implies
$x\leqslant y.$
\end{Proposition}

\begin{proof}
Let $y=s_{i_1}s_{i_2}\cdots s_{i_r}$ and $w_0=t_{j_1}\cdots t_{j_n}$ be reduced expressions.
Since $l(yw_0)=l(y)+l(w_0)$ and $yw_0=s_{i_1}s_{i_2}\cdots s_{i_r}t_{j_1}\cdots t_{j_n}$ is a reduced expression. Then $x\leqslant yw_0$ implies
that $x$ has a reduced expression which is a subexpression of $s_{i_1}s_{i_2}\cdots s_{i_r}t_{j_1}\cdots t_{j_n}.$ Since $l(xw_0)=l(x)+l(w_0)$,
this reduced expression can't end with any $t_{j_p}$. Thus $x\leqslant y.$
\end{proof}

 For $x\in X^+,$ let $V(x)$ be an irreducible rational $G$-module  of
highest weight $x$ and let $S_x$ be the corresponding element defined in \citep{Lusztig1983Singularities}. Then $S_x,x\in X^+$ form an
$\mathcal{A}$-basis of the center of $H.$ For $w\in W_0$ we define
$$E_{d_w}=q^{-\frac{l(d_w)}{2}}\sum\limits_{\mbox{\tiny$\begin{array}{c} y\leqslant d_w\\ l(yw_0)=l(y)+l(w_0)
\end{array}$}}P_{yw_0,d_ww_0}T_y$$
and
$$F_{d_w}=q^{-\frac{l(d_w)}{2}}\sum\limits_{\mbox{\tiny$\begin{array}{c} y\leqslant d_w\\ l(yw_0)=l(y)+l(w_0)
\end{array}$}}P_{yw_0,d_ww_0}T_{y^{-1}}.$$

\begin{Proposition}\label{prop:EF}
\begin{enumerate}[(1)]
    \item $E_{d_w}S_xC_{w_0}F_{d_u}=C_{d_wxw_0d_u^{-1}}$ for any $w,u\in W_0$ and $x\in X^+.$
    \item $S_xS_y=\sum\limits_{z\in X^+}m_{x,y,z}S_z$ for any $x,y\in X^+$, where $m_{x,y,z}$ is defined to be the multiplicity of $V(z)$ in the
tensor product $V(x)\otimes V(y).$
\end{enumerate}
\end{Proposition}

\begin{proof}
(1) See the Corollary 2.11 of \citep{nanhua1990based}.
\\(2) See (8.3) of \citep{Lusztig1983Singularities}.
\end{proof}

As a consequence we have the following lemma.

\begin{Lemma}\label{lem:mu}
For $x\in X^+$ set $x^*=w_0x^{-1}w_0.$\\
(1) For $x,y,z\in X^+,\ m_{x^*,y,z^*}=m_{z,y,x}.$\\
(2) For $y=d_uxw_0,w=d_{u^\prime}x^\prime w_0\in \Gamma_0,$ if $u\neq u^\prime,$ then
$$\mu(y,w)=\delta_{y^{-1},w,w_0}=\sum\limits_{z_1\in X^+}m_{x^*,x^\prime,z_1^*}\delta_{w_0d_u^{-1},d_{u^\prime} w_0,z_1w_0}.$$\\
(3) If $w\neq e\in W_0,$ then $\mu(xw_0,d_ww_0)
=\delta_{w_0,d_ww_0,xw_0}.$
\end{Lemma}

\begin{proof}
(1)
\begin{equation}
\begin{aligned}
m_{x^*,y,z^*}&=\dim\ \text{Hom}_G(V(x^*)\otimes V(y),V(z^*))\\
&=\dim\ \text{Hom}_G(V(x^*)\otimes V(y)\otimes V(z), \mathbb{C})\\
&=\dim\ \text{Hom}_G(V(z)\otimes V(y), V(x))\\
&=m_{z,y,x}.\nonumber
\end{aligned}
\end{equation}
\\(2) See the proof of Theorem 3.1 of \citep{scott2010some}.
\\(3) Note that $c_0$ is the lowest two-sided cell. By (1) and (4) of Proposition \ref{prop:cell} we have that
$$C_{w_0}C_{d_ww_0}=\sum\limits_{t\in X^+}h_{w_0,d_ww_0,tw_0}C_{tw_0}.$$
Multiply $S_{x^*}$ from the left side. By Proposition \ref{prop:EF}(1) we get that
$$C_{x^*w_0}C_{d_ww_0}=h_{w_0,d_ww_0,xw_0}S_{x^*}S_{x}C_{w_0}+\sum\limits_{t\neq x \in X^+}h_{w_0,d_ww_0,tw_0}S_{x^*}S_{t}C_{w_0}.$$
By (1), $m_{x^*,y,e}=m_{e,y,x}=\delta_{x,y}$  for any $y\in X^+.$ Compare the coefficients of $C_{w_0}$ and we have
$$h_{x^*w_0,d_ww_0,w_0}=h_{w_0,d_ww_0,xw_0}.$$
Comparing the coefficients of $q^{\frac{a(w_0)-1}{2}}$ on both sides, we obtain the equality
$$\delta_{x^*w_0,d_ww_0,w_0}=\delta_{w_0,d_ww_0,xw_0}.$$
Then the result follows from (2).
\end{proof}

On the number $m_{x,y,z},$ we have the following useful proposition.

\begin{Proposition}\label{prop:weight}
For $x,y,z\in X^+,$ if $\langle y+\theta,\alpha^{\vee}\rangle\geqslant -1$ for any simple root $\alpha$
and any weight $\theta$ of $V(x)$ then $m_{x,y,z}=\dim\ V(x)_{z-y}.$
\end{Proposition}

\begin{proof}
See Corollary (3.4) of \citep{kumar2010tensor}.
\end{proof}

\section{The case of type $\tilde B_n$ }

In this section we assume that $G=Sp_{2n}(\mathbb{C})$ and $W$ is the corresponding extended affine Weyl group. The Coxeter
diagram of the affine Weyl group is
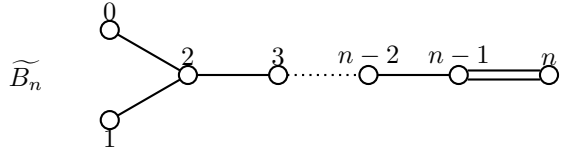
\begin{center}
  \begin{tikzpicture}[scale=.6]
    \draw (-1,0) node[anchor=east]  {$\tilde B_n$ };  
    \draw[thick] (2 cm,0) circle (.2 cm) node [above] {$2$};  
    \draw[xshift=2 cm,thick] (150:2) circle (.2 cm) node [above] {$0$};  
    \draw[xshift=2 cm,thick] (210:2) circle (.2 cm) node [below] {$1$};
    \draw[thick] (4 cm,0) circle (.2 cm) node [above] {$3$};
    \draw[thick] (6 cm,0) circle (.2 cm) node [above] {$n-2$};
    \draw[thick] (8 cm,0) circle (.2 cm) node [above] {$n-1$};
    \draw[thick] (10 cm,0) circle (.2 cm) node [above] {$n$};
    \draw[xshift=2 cm,thick] (150:0.2) -- (150:1.8); 
    \draw[xshift=2 cm,thick] (210:0.2) -- (210:1.8);
    \draw[thick] (2.2,0) --+ (1.6,0);   
    \draw[dotted,thick] (4.2,0) --+ (1.6,0);
    \draw[thick] (6.2,0) --+ (1.6,0);
    \draw[thick] (8.2,0.1) --+ (1.6,0);
    \draw[thick] (8.2,-0.1) --+ (1.6,0);
  \end{tikzpicture}
\end{center}

Note that $X$ is the  weight lattice of the root system $R$. We  use the standard realization of root system (see \citep{Bourbaki1990Lie}or \citep{Humphreys1972Introduction}) to do
calculation . That is, take a standard orthogonal basis $\varepsilon_i\in
E,\ 1\leqslant i\leqslant n,$ such that
$\alpha_1=\varepsilon_1-\varepsilon_2,\alpha_2=\varepsilon_2-\varepsilon_3,\cdots,\alpha_{n-1}=\varepsilon_{n-1}-\varepsilon_{n},
\alpha_n=2\varepsilon_n$ are the simple roots and the root system is
$$R=\{\pm(\varepsilon_{i}\pm\varepsilon_{j}),\pm2\varepsilon_{k}|1\leqslant i\neq j\leqslant n,\ 1\leqslant k\leqslant n\}.$$
Then $\alpha_0=\varepsilon_1+\varepsilon_2.$ Denote the simple reflections by $s_i=s_{\alpha_i},1\leqslant i\leqslant n$ and denote the
fundamental dominant weights by $x_i=x_{\alpha_i},1\leqslant i\leqslant n.$ We have
$x_i=\varepsilon_1+\varepsilon_2+\cdots+\varepsilon_i.$ If we denote the associated inner product of $E$ by $(-,-),$ then $\langle
x,\alpha^{\vee}\rangle=\frac{2(x,\alpha)}{(\alpha,\alpha)}$ for any $x\in X,\ \alpha\in R.$

\begin{Proposition}\label{prop:Bdw}
Keep the notations above. If $n\geqslant 3,$ then for
$w=s_{\varepsilon_1+\varepsilon_3}\in W_0,$ we have $x_2<d_w$ and $l(d_w)=l(x_2)+1.$
\end{Proposition}

\begin{proof}
For $w=s_{\varepsilon_1+\varepsilon_3},$ we have  $d_w=s_{\varepsilon_1+\varepsilon_3}x_1x_3$. Note that $s_2s_0s_2=d_wx_2^{-1}.$ Thus $d_w=s_2s_0s_2x_2.$ By the length formula, we have $l(d_w)=4n-3$, and $l(x_2)=\langle x_2,2\rho^{\vee}\rangle=(\varepsilon_1+\varepsilon_2,(2n-1)\varepsilon_1+(2n-3)\varepsilon_2+\cdots+\varepsilon_n)=4n-4.$ Then $l(d_w)=l(x_2)+1.$

It is easy to see $s_2x_2<x_2$ by computing their lengths. Thus $x_2$ has a reduced expression which is a subexpression of a reduced
expression of $d_w=s_2s_0s_2x_2.$ Hence $x_2<d_w.$
\end{proof}

Now we can state the main result in this section.

\begin{Theorem}\label{thm:Bmu}
Let $G=Sp_{2n}(\mathbb{C}).$ If $n\geqslant 4,$ then for
$x=\prod\limits_{i=1}^nx_i^{a_i},a_i\geqslant1-\delta_{i,n}$ and $w=s_{\varepsilon_1+\varepsilon_3},$ we have
$$\mu(xw_0,d_wxw_0)=n.$$
\end{Theorem}

\begin{proof}
By Lemma \ref{lem:mu} (2) we have
$$\mu(xw_0,d_wxw_0)=\sum\limits_{z_1\in X^+}m_{x^*,x,z_1^*}\delta_{w_0,d_{w} w_0,z_1w_0}.$$

Assume that $\delta_{w_0,d_ww_0,z_1w_0}\neq 0.$ Note that $w_0,\ d_ww_0\in W_a,$ so $z_1\in X^+\cap \mathbb{Z}R.$ Since $w_0\not\sim_L (d_ww_0)^{-1}$, so $\gamma_{w_0,d_ww_0,z_1w_0}= 0.$ By Proposition \ref{prop:length} we have
$a(w_0)-1\leqslant l(w_0)+l(d_ww_0)-l(z_1w_0).$ Thus $l(z_1)\leqslant l(d_w)+1=4n-2.$ Since $l(x_i)=2ni-i^2$, we have $l(x_1^2)=4n-2,\ l(x_2)=4n-4,\
l(x_i)> 4n-2 $ for $3\leqslant i\leqslant n,$ and $l(x_1x_2)=6n-5>4n-2.$ (Here we need $n\geqslant 4$.) Therefore $z_1\in
\{e,x_1^2,x_2\}$.

By Lemma \ref{lem:mu} (3) and Proposition \ref{prop:Bdw} we have
$$\delta_{w_0,d_ww_0,x_2w_0}=\mu(x_2w_0,d_ww_0)=1.$$
Similarly we have
$$\delta_{w_0,d_ww_0,x_1^2w_0}=\mu(x_1^2w_0,d_ww_0)=0,$$
since $d_w\not < x_1^2$ and $x_1^2\not < d_w.$ In fact,
since $l(d_w)=l(x_1^2)-1,\ s_1x_1^2<x_1^2$ and $s_1d_w>d_w,$ so if $d_w<x_1^2,$ we must have
$d_w=s_1x_1^2,$
which is impossible by an easy computation.
Hence
$$\mu(xw_0,d_wxw_0)=m_{x^*,x,x_2^*}+m_{x^*,x,e}\delta_{w_0,d_ww_0,w_0}=m_{x_2,x,x}+\mu(w_0,d_ww_0).$$

Now we compute $m_{x_2,x,x}.$ Obviously, the set of weights of $V(x_2)$ is $\{u(x_2),0|u\in W_0\}=\{\pm(\varepsilon_i\pm\varepsilon_j),0|i\neq
j\}.$ For any $u\in W_0,$ $\langle u(x_2),\alpha_k^{\vee}\rangle\geqslant -2$ for all $1\leqslant k< n$ and $\langle
u(x_2),\alpha_n^{\vee}\rangle\geqslant -1.$ Thus by Proposition \ref{prop:weight}, $m_{x_2,x,x}=\dim\ V(x_2)_0.$
Using Freudenthal's weight multiplicity formula we obtain
$\dim\ V(x_2)_0=n-1.$

We are now reduced to show that $\mu(w_0,d_ww_0)=1$, which is  Lemma \ref{lem:B3} to be proved next. The theorem is proved.
\end{proof}

Now we are going to prove $\mu(w_0,d_ww_0)=1.$ First we need to work out a reduced expression of $d_w.$ Let
$\omega$ be the unique nontrivial element in $\Omega$. Then $\omega^2=1,\ \omega s_1 \omega=s_0,\ \omega s_i
\omega=s_i,\ 2\leqslant i\leqslant n.$ Since $x_1=\varepsilon_1,\ x_2=\varepsilon_1+\varepsilon_2,$ we have $x_1=\omega
s_0s_2s_3...s_n...s_3s_2s_0,$ thus
\begin{equation}
\begin{aligned}
x_2&=x_1s_1x_1s_1\\
&=s_2s_3\cdots s_n\ s_1s_2\cdots s_{n-2}s_{n-1}s_{n-2}\cdots s_2s_1\ s_n \cdots s_3s_2\ s_0.\nonumber
\end{aligned}
\end{equation}
We then  get a reduced expression of $d_w$ as follows
$$d_w=s_2s_0\ s_3\cdots s_n\ s_1s_2\cdots s_{n-2}s_{n-1}s_{n-2}\cdots s_2s_1\ s_n \cdots s_3s_2\ s_0.$$

Set $v=d_w,v_1=s_2d_w$ and $v_2=s_0s_2d_w.$ Since $w_0,vw_0\in \mathscr{D}_L(s_0,s_2)=\{w\in W_a|\#L(w)\cap \{s_0,s_2\}=1\},$ using the star
operator we see $\mu(w_0,vw_0)=\mu(s_0w_0,v_1w_0).$

Using Proposition \ref{prop:kl} (1) we get
\begin{equation}
P_{s_0w_0,v_1w_0}=P_{w_0,v_2w_0}+qP_{s_0w_0,v_2w_0}-\sum\limits_{\mbox{\tiny$\begin{array}{c}s_0w_0\leqslant z\prec v_2w_0\\
s_0z<z\end{array}$}}\mu(z,v_2w_0)q^{\frac{1}{2}(l(v_1w_0)-l(z))}P_{s_0w_0,z}.
\tag{$*$}
\end{equation}

We first show that the summation in formula $(*)$ is empty. That is, there exists no $z\in W$ such that $s_0w_0\leqslant z\prec v_2w_0$ and
$s_0z<z$. Assume that $z\in W$ satisfies the conditions. Thanks to $s_0v_2w_0> v_2w_0,$ the conditions $z\prec v_2w_0$ and $s_0z<z$ imply
that $z\leqslant_L v_2w_0.$ Since $v_2w_0\in \Gamma_0$, we have $z\in \Gamma_0.$ In particular, $z=uw_0$ for some $u\in W$ and
$l(z)=l(u)+l(w_0).$ Since $s_0w_0\leqslant z\leqslant v_2w_0,$ and
$$v_2=s_3\cdots s_n\ s_1s_2\cdots s_{n-2}s_{n-1}s_{n-2}\cdots s_2s_1\ s_n \cdots s_3s_2\ s_0,$$
by Proposition \ref{prop:leftcell} we see that $z$ must be one of the following elements:
\begin{equation}
\begin{aligned}
&s_0w_0;\\
&m_i=s_i\cdots s_3s_2s_0w_0,\ 2\leqslant i\leqslant n;\\
&m_{j,i}=s_j\cdots s_2s_1\ \ s_i\cdots s_2s_0w_0,\ 2\leqslant i\leqslant n, 1\leqslant j<i;\\
&g_k=s_k\cdots s_{n-1}\ \ s_n\cdots s_2s_0w_0,\ 1\leqslant k\leqslant n-1;\\
&w_{k,j}=s_k\cdots s_{n-1}s_j\cdots s_2s_1\ \ s_n\cdots s_2s_0w_0,\ 1\leqslant k\leqslant n-1, 1\leqslant j\leqslant n-2;\\
&z_{t,k}=s_t\cdots s_n\ \ s_k\cdots s_{n-1}\cdots s_1\ \ s_n\cdots s_2s_0w_0,\ 3\leqslant t\leqslant n,1\leqslant k<t.
\nonumber
\end{aligned}
\end{equation}

Since $\mathcal{L}(v_2w_0)=S-\{s_2,s_0\}$, combining the condition $s_0z<z$ and Proposition \ref{prop:kl} (4), we see that $z$ is one of the following
elements:
 $$s_0w_0,\quad w_{3,1}.$$

Since $l(v_2w_0)-l(s_0w_0)=4n-6$ is even, so $\mu(s_0w_0,v_2w_0)=0.$ Note that $\mathcal{L}(w_{3,1})=S-\{s_2\}.$ Using the star operator with
respect to the pair $\{s_2,s_3\}$, we see that
$$\mu(w_{3,1},v_2w_0)=\mu(w_{2,1},z_{4,1})=0.$$
Therefore we have
\begin{equation}
P_{s_0w_0,v_1w_0}=P_{w_0,v_2w_0}+qP_{s_0w_0,v_2w_0}.
\tag{\#}
\end{equation}



\begin{Lemma}\label{lem:B1}
 Let $v_2$ be as above. If $n\geqslant 4$ then
$$P_{w_0,v_2w_0}=q^{2n-3}+q^{2n-4}+\cdots+q+1.$$
\end{Lemma}

\begin{proof}
If $n=4,5,$ the calculation is easy and we omit it. When $n\geqslant 6$, the case that $n$ is even and the case that $n$
is odd are similar. In view of this, we only consider the case that $n\geqslant 6$ is even. Recall that we have a reduced expression
$$v_2=s_3\cdots s_n\ s_1s_2\cdots s_{n-2}s_{n-1}s_{n-2}\cdots s_2s_1\ s_n \cdots s_3s_2\ s_0.$$
In the following we shall use an analysis similar to that for $(*)$ to compute $P_{w_0,v_2w_0}.$ Starting with the above reduced expression, every time we cancel the leftmost simple reflection. Thus we obtain a series of equalities (like ($\#$)) about Kazhdan-Lusztig polynomials and some $\mu(y,w)$. Then we apply flexibly the properties in Proposition \ref{prop:kl} to compute what we need to know in these equalities.  Note that $v_2w_0=z_{3,1}.$

{\bf (1)} Recall that for $4\leqslant i\leqslant n,$
$$z_{i,1}=s_i\cdots s_n\ s_1s_2\cdots s_{n-2}s_{n-1}s_{n-2}\cdots s_2s_1\ s_n \cdots s_3s_2\ s_0w_0,$$
and then $z_{i,1}=s_{i-1}z_{i-1,1}.$

For $4\leqslant i\leqslant n,$ since $s_{i-1}w_0=w_0s_{i-1}\leqslant w_0$ and $z_{i,1}s_{i-1}<z_{i,1}$, by Proposition \ref{prop:kl} (2) we have
$$P_{s_{i-1}w_0,z_{i,1}}=P_{w_0s_{i-1},z_{i,1}}=P_{w_0,z_{i,1}}.$$

Since $z_{i,1}=s_{i-1}z_{i-1,1}\leqslant z_{i-1,1}$, using the recursive formula in Proposition \ref{prop:kl} (1) and through an analysis similar to that for $(*)$, noting that $\mathcal{L}(z_{i,1})=S-\{s_{i-1},s_0\},l(z_{i,1})=4n-2-i+l(w_0)$,  we get
\begin{equation}
P_{w_0,z_{i-1,1}}=
\begin{cases}
(1+q)P_{w_0,z_{i,1}}-qP_{w_0,z_{i+1,1}},&i=4,6,\cdots,n-2,\\
(1+q)P_{w_0,z_{i,1}}-\mu(w_0,z_{i,1})q^{2n-1-\frac{i-1}{2}}-qP_{w_0,z_{i+1,1}}, &i=5,7,\cdots,n-1,\\
(1+q)P_{w_0,z_{i,1}}-qP_{w_0,w_{1,n-2}},&i=n,
\end{cases}
\tag{3.3.1}
\end{equation}
where $w_{1,n-2}=s_1s_2\cdots s_{n-2}s_{n-1}s_{n-2}\cdots s_2s_1\ s_n \cdots s_3s_2\ s_0w_0$.

{\bf (2)} Similarly, since $w_{1,n-2}=s_nz_{n,1}\leqslant z_{n,1}$, using the recursive formula in Proposition \ref{prop:kl} (1) and through  an analysis similar to that for $(*)$, noting that $\mathcal{L}(w_{1,n-2})=S-\{s_n,s_0\},\ l(w_{1,n-2})=3n-3+l(w_0)$,  we get
\begin{equation}
P_{w_0,z_{n,1}}=(1+q)P_{w_0,w_{1,n-2}}-\mu(w_0,w_{1,n-2})q^{\frac{3n}{2}-1}-qP_{w_0,w_{1,n-3}},
\tag{3.3.2}
\end{equation}
where $w_{1,n-3}=s_1s_2\cdots s_{n-2}s_{n-1}s_{n-3}\cdots s_2s_1\ s_n \cdots s_3s_2\ s_0w_0.$

Recall that for $2\leqslant i\leqslant n-2$,
$$w_{i,n-2}=s_i\cdots s_{n-2}s_{n-1}s_{n-2}\cdots s_2s_1\ s_n \cdots s_3s_2\ s_0w_0,$$
and then $w_{i,n-2}=s_{i-1}w_{i-1,n-2}\leqslant w_{i-1,n-2}.$

Using the recursive formula in Proposition \ref{prop:kl} (1) and through  an analysis similar to that for $(*)$, noting that $\mathcal{L}(w_{i,n-2})=S-\{s_n,s_{i-1},s_0\},\ l(w_{i,n-2})=3n-i-2+l(w_0)$, we get
\begin{equation}
\begin{aligned}
&P_{w_0,w_{i-1,n-2}}=\\
&\begin{cases}
(1+q)P_{w_0,w_{i,n-2}}-\mu(m_{n-2,n-1},w_{i,n-2})q^{\frac{n-i+2}{2}}P_{w_0,m_{n-2,n-1}}-qP_{w_0,w_{i+1,n-2}}, &i=2,4,\cdots,n-4,\\
(1+q)P_{w_0,w_{i,n-2}}-\mu(w_0,w_{i,n-2})q^{\frac{3n-i+1}{2}-1}-qP_{w_0,w_{i+1,n-2}},&i=3,5,\cdots,n-3,\\
(1+q)P_{w_0,w_{i,n-2}}-\mu(m_{n-2,n-1},w_{i,n-2})q^{\frac{n-i+2}{2}}P_{w_0,m_{n-2,n-1}}-qP_{w_0,m_{n-1,n}},&i=n-2,
\end{cases}
\end{aligned}
\tag{3.3.3}
\end{equation}
where $m_{n-2,n-1}=s_{n-2}\cdots s_1\ s_{n-1}\cdots s_2\ s_0w_0,$ and
$m_{n-1,n}=s_{n-1}\cdots s_1\ s_{n}\cdots s_2\ s_0w_0.$

{\bf (3)} Since $m_{n-1,n}=s_{n-2}w_{n-2,n-2}\leqslant w_{n-2,n-2}$, using the recursive formula in Proposition \ref{prop:kl} (1) and through  an analysis similar to that for $(*)$, noting that $\mathcal{L}(m_{n-1,n})=S-\{s_n,s_{n-2},s_0\},\ l(m_{n-1,n})=2n-1+l(w_0)$, we get
\begin{equation}
P_{w_0,w_{n-2,n-2}}=(1+q)P_{w_0,m_{n-1,n}}-\mu(w_0,m_{n-1,n})q^{n}-qP_{w_0,m_{n-2,n}},
\tag{3.3.4}
\end{equation}
where $m_{n-2,n}=s_{n-2}\cdots s_1\ s_n\cdots s_2\ s_0w_0.$

Since $m_{n-2,n}=s_{n-1}m_{n-1,n}\leqslant m_{n-1,n}$, using the recursive formula in Proposition \ref{prop:kl} (1) and through  an analysis similar to that for $(*)$, noting that $\mathcal{L}(m_{n-2,n})=S-\{s_{n-1},s_0\},\ l(m_{n-2,n})=2n-2+l(w_0)$, we get
\begin{equation}
P_{w_0,m_{n-1,n}}=(1+q)P_{w_0,m_{n-2,n}}-qP_{w_0,m_{n-2,n-1}}.
\tag{3.3.5}
\end{equation}

Recall that for $3\leqslant i\leqslant n-2,$
$$m_{n-i,n}=s_{n-i}\cdots s_1\ s_{n}\cdots s_2\ s_0w_0,$$
and then $m_{n-i,n}=s_{n-i+1}m_{n-i+1,n}\leqslant m_{n-i+1,n}.$

using the recursive formula in Proposition \ref{prop:kl} (1) and through  an analysis similar to that for $(*)$, noting that $\mathcal{L}(m_{n-i,n})=S-\{s_{n-1},s_{n-i+1},s_0\},\ l(m_{n-i,n})=2n-i+l(w_0)$, we get
\begin{equation}
P_{w_0,m_{n-i+1,n}}=
\begin{cases}
(1+q)P_{w_0,m_{n-i,n}}-\mu(w_0,m_{n-i,n})q^{n-\frac{i-1}{2}}-qP_{w_0,m_{n-i-1,n}}, &i=3,5,\cdots,n-3,\\
(1+q)P_{w_0,m_{n-i,n}}-qP_{w_0,m_{n-i-1,n}},&i=4,6,\cdots,n-2,
\end{cases}
\tag{3.3.6}
\end{equation}
where $m_{1,n}=s_1\ s_n\cdots s_2\ s_0w_0.$

Since $m_{1,n}=s_2m_{2,n}\leqslant m_{2,n}$, using the recursive formula in Proposition \ref{prop:kl} (1) and through  an analysis similar to that for $(*)$, noting that $\mathcal{L}(m_{1,n})=S-\{s_{n-1},s_2\},\ l(m_{1,n})=n+1+l(w_0),$ we get
\begin{equation}
P_{w_0,m_{2,n}}=(1+q)P_{w_0,m_{1,n}}-qP_{w_0,m_n},
\tag{3.3.7}
\end{equation}
where $m_{n}=s_n\cdots s_2\ s_0w_0.$

{\bf (4)} Since $m_{n}=s_1m_{1,n}\leqslant m_{1,n}$, using the recursive formula in Proposition \ref{prop:kl} (1) and through  an analysis similar to that for $(*)$, noting that $\mathcal{L}(m_n)=S-\{s_1,s_{n-1}\},l(m_n)=n+l(w_0),$ we get
\begin{equation}
P_{w_0,m_{1,n}}=(1+q)P_{w_0,m_{n}}.
\tag{3.3.8}
\end{equation}

{\bf (5)} Using the recursive formula in Proposition \ref{prop:kl} (1) and through  an analysis similar to that for $(*)$,  we can easily see  $P_{w_0,m_{n}}=1.$  Combining formulas (3.3.6)-(3.3.8), by a simple computation we get
\begin{equation}
P_{w_0,m_{n-2,n}}=q^{n-2}+q^{n-3}+\cdots+q+1.
\tag{3.3.9}
\end{equation}
To complete the calculations in equalities (3.3.1)-(3.3.5), it remains to figure out the followings:
\begin{align*}
&\mu(m_{n-2,n-1},w_{i,n-2}),i=2,4,\cdots,n-2;\\
&P_{w_0,m_{n-2,n-1}},\quad P_{w_0,w_{1,n-3}}.
\end{align*}

{\bf (6)} Recall that for $2\leqslant i\leqslant n-2$,
$$w_{i,n-2}=s_i\cdots s_{n-2}s_{n-1}s_{n-2}\cdots s_2s_1\ s_n \cdots s_3s_2\ s_0w_0,$$
$$w_{n-2,n-3}=s_{n-2}s_{n-1}s_{n-3}\cdots s_1\ s_n\cdots s_2\ s_0w_0,$$
and
$$m_{n-2,n-1}=s_{n-2}\cdots s_1\ s_{n-1}\cdots s_2\ s_0w_0.$$

Since $s_{n-1}m_{n-2,n-1}=m_{n-2,n-1}s_1$ and $ w_{n-2,n-3}s_1\leqslant w_{n-2,n-3}$,  we have
$$P_{s_{n-1}m_{n-2,n-1},w_{n-2,n-3}}=P_{m_{n-2,n-1},w_{n-2,n-3}}.$$

Since $w_{n-2,n-3}=s_{n-1}w_{n-2,n-2}\leqslant w_{n-2,n-2}$,  using the recursive formula in Proposition \ref{prop:kl} (1) and through  an analysis similar to that for $(*)$, noting that $\mathcal{L}(w_{n-2,n-3})=S-\{s_{n-1},s_{n-3},s_0\},\ l(w_{n-2,n-3})=2n-1+l(w_0)$, we get
$$P_{m_{n-2,n-1},w_{n-2,n-2}}=(1+q)P_{m_{n-2,n-1},w_{n-2,n-3}}.$$
But $l(w_{n-2,n-3})-l(m_{n-2,n-1})=2$ implies $P_{m_{n-2,n-1},w_{n-2,n-3}}=1,$ thus
\begin{equation}
P_{m_{n-2,n-1},w_{n-2,n-2}}=q+1.
\tag{3.3.10}
\end{equation}
Therefore
\begin{equation}
\mu(m_{n-2,n-1},w_{n-2,n-2})=1.
\tag{3.3.11}
\end{equation}

Since $s_{i}m_{n-2,n-1}=m_{n-2,n-1}s_{i+2}$ for $1\leqslant i\leqslant n-3$, from (3.3.3) we see that
\begin{equation}
\begin{aligned}
&P_{m_{n-2,n-1},w_{i-1,n-2}}=\\
&\begin{cases}
(1+q)P_{m_{n-2,n-1},w_{i,n-2}}-\mu(m_{n-2,n-1},w_{i,n-2})q^{\frac{n-i+2}{2}}-qP_{m_{n-2,n-1},w_{i+1,n-2}}, &i=2,4,\cdots,n-4,\\
(1+q)P_{m_{n-2,n-1},w_{i,n-2}}-qP_{m_{n-2,n-1},w_{i+1,n-2}},&i=3,5,\cdots,n-3,\\
(1+q)P_{m_{n-2,n-1},w_{i,n-2}}-\mu(m_{n-2,n-1},w_{i,n-2})q^{\frac{n-i+2}{2}}-qP_{m_{n-2,n-1},m_{n-1,n}}, &i=n-2.
\end{cases}
\nonumber
\end{aligned}
\end{equation}
But $l(m_{n-1,n})-l(m_{n-2,n-1})=2$ implies $P_{m_{n-2,n-1},m_{n-1,n}}=1,$ thus for $2\leqslant i\leqslant n-4$,
$$P_{m_{n-2,n-1},w_{i,n-2}}=q+1,$$
and
\begin{equation}
\mu(m_{n-2,n-1},w_{i,n-2})=0.
\tag{3.3.12}
\end{equation}

{\bf (7)} Recall that for $3\leqslant i\leqslant n-1$,
$$m_{n-i,n-1}=s_{n-i}\cdots s_1\ s_{n-1}\cdots s_2\ s_0w_0.$$
Thus $m_{n-i,n-1}=s_{n-i+1}m_{n-i+1,n-1}\le m_{n-i+1,n-1}.$ Obviously $l(m_{n-i,n-1})=2n-i-1+l(w_0)$ and
\[
\mathcal{L}(m_{n-i,n-1})=
\begin{cases}
S-\{s_n,s_{n-i+1},s_0\}, &3\leqslant i\leqslant n-2\\
S-\{s_n,s_2\}, &i=n-1.
\end{cases}
\]
In addition, for $m_{n-1}=s_{n-1}\cdots s_2\ s_0w_0,$ we have
$$\mathcal{L}(m_{n-1})=S-\{s_1,s_{n}\},\ l(m_{n-1})=n-1+l(w_0).$$
 using the recursive formula in Proposition \ref{prop:kl} (1) and through  an analysis similar to that for $(*)$,   we see that
\begin{align*}
&P_{w_0,m_{n-i+1,n-1}}=\\
&\begin{cases}
(1+q)P_{w_0,m_{n-i,n-1}}-q(P_{w_0,m_{n-i,n-2}}+P_{w_0,m_{n-i-1,n-1}}),&i=3,\\
(1+q)P_{w_0,m_{n-i,n-1}}-\mu(w_0,m_{n-i-1,n-1})q^{n-\frac{i}{2}}-qP_{w_0,m_{n-i-1,n-2}},&i=4,6,\cdots,n-2,\\
(1+q)P_{w_0,m_{n-i,n-1}}-qP_{w_0,m_{n-i-1,n-2}},&i=5,7,\cdots,n-3,\\
(1+q)P_{w_0,m_{n-i,n-1}}-qP_{w_0,m_{n-1}},&i=n-1,\\
(1+q)P_{w_0,m_{n-1}},&i=n.
\end{cases}
\end{align*}

An easy computation leads to $P_{w_0,m_{n-1}}=1.$ By a direct calculation, we see that
\begin{equation}
P_{w_0,m_{n-i,n-1}}=q^{n-i}+q^{n-i-1}+\cdots+q+1,\ 3\leqslant i\leqslant n-1.
\tag{3.3.13}
\end{equation}
Thus
$$P_{w_0,m_{n-2,n-1}}=(1+q)(q^{n-3}+\cdots+1)-q(P_{w_0,m_{n-3,n-2}}+(q^{n-4}+\cdots+1)).$$
 Similarly we get
$$P_{w_0,m_{n-3,n-2}}=(1+q)(q^{n-4}+\cdots+1)-q(P_{w_0,m_{n-4,n-3}}+(q^{n-5}+\cdots+1)).$$
By easy calculations we get  $P_{w_0,m_{2,3}}=q^2+1$ and $P_{w_0,m_{3,4}}=q^2+1.$ Using induction on $n\geqslant 6$ even, we can prove that:
\begin{equation}
P_{w_0,m_{n-2,n-1}}=q^{n-2}+q^{n-4}+\cdots+q^2+1,
\tag{3.3.14}
\end{equation}
and
\begin{equation}
P_{w_0,m_{n-3,n-2}}=q^{n-4}+q^{n-6}+\cdots+q^2+1.
\tag{3.3.15}
\end{equation}

{\bf (8)} Recall that for $2\leqslant i\leqslant n-1,$
$$w_{i,n-3}=s_{i}\cdots s_{n-2}s_{n-1}s_{n-3}\cdots s_1\ s_n\cdots s_2\ s_0w_0.$$
Thus $w_{i,n-3}=s_{i-1}w_{i-1,n-3}\le w_{i-1,n-3}.$ Obviously $l(w_{i,n-3})=3n-i-3+l(w_0)$ and
\[
\mathcal{L}(w_{i,n-3})=
\begin{cases}
S-\{s_{n-1},s_{i-1},s_0\} &2\leqslant i\leqslant n-2\\
S-\{s_{n-2},s_0\} &i=n-1.
\end{cases}
\]
In addition, for $m_{n-3,n}=s_{n-3}\cdots s_1\ s_{n}\cdots s_2\ s_0w_0=s_{n-1}w_{n-1,n-3}$, we have
$$\mathcal{L}(m_{n-3,n})=S-\{s_{n-1},s_{n-2},s_0\},\ l(m_{n-3,n})=2n-3+l(w_0).$$
Using the recursive formula in Proposition \ref{prop:kl} (1) and through  an analysis similar to that for $(*)$, we get
\begin{equation}
\begin{aligned}
&P_{w_0,w_{i-1,n-3}}=\\
&\begin{cases}
(1+q)P_{w_0,w_{i,n-3}}-\mu(w_0,w_{i,n-3})q^{\frac{3n-2-i}{2}}-\mu(m_{n-2,n},w_{i,n-3})q^{\frac{n-i}{2}}P_{w_0,m_{n-2,n}}-qP_{w_0,w_{i+1,n-3}},\\
\ \ \ \ \ \ \ \ \ \ \ \ \ \ \ \ \ \ \ \ \ \ \ \ \ \ \ \ \ \ \ \ \ \ \ \ \ \ \ \ \ \ \ \ \ \ \ \ \ \ \ \ \ \ \ \ \ \ \ \ \ \ \ \ \ \ \ \ \ \ \ \
 \ \ \ \ i=2,4,\cdots,n-2,\\
(1+q)P_{w_0,w_{i,n-3}}-\mu(m_{n-3,n-2},w_{i,n-3})q^{\frac{n-i+3}{2}}P_{w_0,m_{n-3,n-2}}-qP_{w_0,w_{i+1,n-3}},\\
\ \ \ \ \ \ \ \ \ \ \ \ \ \ \ \ \ \ \ \ \ \ \ \ \ \ \ \ \ \ \ \ \ \ \ \ \ \ \ \ \ \ \ \ \ \ \ \ \ \ \ \ \ \ \ \ \ \ \ \ \ \ \ \ \ \ \ \ \ \ \ \
 \ \ \ \
i=3,5,\cdots,n-3,\\
(1+q)P_{w_0,w_{i,n-3}}, \ \ \ \ \ \ \ \ \ \ \ \ \ \ \ \ \ \ \ \ \ \ \ \ \ \ \ \ \ \ \ \ \ \ \ \ \ \ \ \ \ \ \ \ \ \ \ \ \ \ \ \ \ \ i=n-1,\\
(1+q)P_{w_0,m_{n-3,n}}-\mu(w_0,m_{n-3,n})q^{n-1}-qP_{w_0,m_{n-3,n-1}},\ \ \ \ \ \ i=n.
\end{cases}
\nonumber
\end{aligned}
\end{equation}

For $m_{n-3,n-1}=s_{n-3}\cdots s_1\ s_{n-1}\cdots s_2\ s_0w_0=s_nm_{n-3,n}\leqslant m_{n-3,n},$ we have
$$\mathcal{L}(m_{n-3,n-1})=S-\{s_n,s_{n-2},s_0\},l(m_{n-3,n-1})=2n-4+l(w_0).$$
Using the recursive formula in Proposition \ref{prop:kl} (1) and through  an analysis similar to that for $(*)$, we get
\[
P_{w_0,m_{n-3,n}}=(1+q)P_{w_0,m_{n-3,n-1}}-\mu(m_{n-4,n-3},m_{n-3,n-1})q^2P_{w_0,m_{n-4,n-3}}-qP_{w_0,m_{n-3,n-2}}.
\]

Then from (3.3.13),(3.3.9) and (3.3.14)-(3.3.15), we see respectively that
$$P_{w_0,m_{n-3,n-1}}=q^{n-2}+q^{n-3}+\cdots+q+1;$$
$$P_{w_0,m_{n-2,n}}=q^{n-2}+q^{n-3}+\cdots+q+1;$$
and
$$P_{w_0,m_{n-3,n-2}}=q^{n-4}+q^{n-6}+\cdots+q^2+1,$$
$$P_{w_0,m_{n-4,n-3}}=q^{n-4}+q^{n-6}+\cdots+q^2+1.$$

Now we use star operators to compute the $\mu(m_{i,j},w_{k,l})$ occurring in the above equalities. First we have
$$\mu(m_{n-2,n},w_{n-2,n-3})=1,$$
since $l(w_{n-2,n-3})-l(m_{n-2,n})=1.$
Using the star operator with respect to the pair $\{s_{n-2},s_{n-1}\}$, we see that for $i=2,4,\cdots,n-4$,
$$\mu(m_{n-2,n},w_{i,n-3})=\mu(m_{n-1,n},w_{i,n-4}).$$
\begin{itemize}
  \item If $i=2,4,\cdots,n-6$, since $s_n\in \mathcal{L}(w_{i+1,n-4})-\mathcal{L}(m_{n-1,n})$ and $l(w_{i+1,n-4})-l(m_{n-1,n})>1$, by Proposition \ref{prop:kl} (4)
we see that
$$\mu(m_{n-2,n},w_{i,n-3})=0.$$
  \item If $i=n-4,$ since $m_{n-1,n}\not\leqslant w_{n-4,n-4},$ we see that
$$\mu(m_{n-2,n},w_{n-4,n-3})=0.$$

\end{itemize}

Similarly, using the same star operator we can obtain that
\[
\mu(m_{n-3,n-2},w_{i,n-3})=\mu(m_{n-3,n-1},w_{i,n-4})=0, i=3,5,\cdots,n-3
\]
and
$$\mu(m_{n-4,n-3},m_{n-3,n-1})=\mu(m_{n-4,n-2},m_{n-3,n-2})=1.$$

Combining the previous paragraphs, by a direct computation we get
\begin{equation}
P_{w_0,w_{1,n-3}}=q^{n-2}+q^{n-3}+\cdots+q+1.
\tag{3.3.16}
\end{equation}

{\bf (9)} Finally using the identities  (3.3.11),(3.3.12),(3.3.14),(3.3.16) we complete the calculations in (3.3.1)-(3.3.5), and prove the formula in the lemma.
\end{proof}

\begin{Lemma}\label{lem:B2}
Let $v_2$ be as in Lemma \ref{lem:B1}. If $n\geqslant 4$ then
$$\deg P_{s_0w_0,v_2w_0}<2n-4.$$
\end{Lemma}

\begin{proof}
This proof is similar to the one for Lemma \ref{lem:B1}. As the proof of Lemma \ref{lem:B1}, we shall only consider the case that
$n\geqslant 6$ is even.  Using the recursive formula in Proposition \ref{prop:kl} (1) and through  an analysis similar to that for $(*)$, especially noting that
$v_2w_0=z_{3,1}, w_{3,n-2}=s_2w_{2,n-2}$ and $m_{1,n}=s_2m_{2,n}$, we obtain the following formulas (1)-(4):

{\bf (1)}
\begin{equation}
P_{s_0w_0,z_{i-1,1}}=
\begin{cases}
(1+q)P_{s_0w_0,z_{i,1}}-qP_{s_0w_0,z_{i+1,1}},&4\leqslant i\leqslant n-1,\\
(1+q)P_{s_0w_0,z_{i,1}}-qP_{s_0w_0,w_{1,n-2}},&i=n.
\end{cases}
\tag{3.4.1}
\end{equation}

{\bf (2)}
\begin{equation}
P_{s_0w_0,z_{n,1}}=(1+q)P_{s_0w_0,w_{1,n-2}}-qP_{s_0w_0,w_{1,n-3}}.
\tag{3.4.2}
\end{equation}

\begin{equation}
\begin{aligned}
&P_{w_0,w_{i-1,n-2}}=\\
&\begin{cases}
(1+q)P_{s_0w_0,w_{i,n-2}}-\mu(m_{n-2,n-1},w_{i,n-2})q^{\frac{n-i+2}{2}}P_{s_0w_0,m_{n-2,n-1}}-qP_{s_0w_0,w_{i+1,n-2}},\\
\ \ \ \ \ \ \ \ \ \ \ \ \ \ \ \ \ \ \ \ \ \ \ \ \ \ \ \ \ \ \ \ \ \ \ \ \ \ \ \ \ \ \ \ \ \ \ \ \ \ \ \ \ \ \ \ \ \ \ \ \ \ \ \ \ \ \ \ \ \ \ \ \ \ \ \ \
i=2,4,\cdots,n-4,\\
P_{s_0w_0,w_{i,n-2}}+qP_{s_2s_0w_0,w_{i,n-2}}-qP_{w_0,w_{i+1,n-2}},
\ \ \ \ \ \ \ \ \ \ \ \ \ \ \ \ \ \ \ \ \ i=3,\\
(1+q)P_{s_0w_0,w_{i,n-2}}-qP_{s_0w_0,w_{i+1,n-2}},
\ \ \ \ \ \ \ \ \ \ \ \ \ \ \ \ \ \ \ \ \ \ \ \ \ \ \ \ \ \ \ i=5,7,\cdots,n-3,\\
(1+q)P_{s_0w_0,w_{i,n-2}}-\mu(m_{n-2,n-1},w_{i,n-2})q^{\frac{n-i+2}{2}}P_{s_0w_0,m_{n-2,n-1}}-qP_{s_0w_0,m_{n-1,n}}, \\
\ \ \ \ \ \ \ \ \ \ \ \ \ \ \ \ \ \ \ \ \ \ \ \ \ \ \ \ \ \ \ \ \ \ \ \ \ \ \ \ \ \ \ \ \ \ \ \ \ \ \ \ \ \ \ \ \ \ \ \ \ \ \ \ \ \ \ \ \ \ \ \ \ \ \ \ \ i=n-2.
\end{cases}
\end{aligned}
\tag{3.4.3}
\end{equation}

{\bf (3)}
\begin{equation}
P_{s_0w_0,w_{n-2,n-2}}=(1+q)P_{s_0w_0,m_{n-1,n}}-qP_{s_0w_0,m_{n-2,n}}.
\tag{3.4.4}
\end{equation}

\begin{equation}
P_{s_0w_0,m_{n-1,n}}=(1+q)P_{s_0w_0,m_{n-2,n}}-qP_{s_0w_0,m_{n-2,n-1}}.
\tag{3.4.5}
\end{equation}

\begin{equation}
P_{s_0w_0,m_{n-i+1,n}}=
(1+q)P_{s_0w_0,m_{n-i,n}}-qP_{s_0w_0,m_{n-i-1,n}}, 3\leqslant i\leqslant n-2.
\tag{3.4.6}
\end{equation}

{\bf (4)}
\begin{equation}
P_{s_0w_0,m_{2,n}}=P_{s_0w_0,m_{1,n}}+qP_{s_2s_0w_0,m_{1,n}}-qP_{s_0w_0,m_{n}}.
\tag{3.4.7}
\end{equation}

\begin{equation}
P_{s_0w_0,m_{1,n}}=(1+q)P_{s_0w_0,m_{n}}.
\tag{3.4.8}
\end{equation}

{\bf (5)}
Since for $m_n=s_n\cdots s_2s_0w_0,$ $\mathcal{L}(m_n)=S-\{s_{n-1},s_1\}$, by Proposition \ref{prop:kl} (2) we see that
$$P_{s_0w_0,m_{n}}=P_{m_{n-2},m_{n}}=1.$$
Thus
$$P_{s_0w_0,m_{1,n}}=q+1.$$

Similarly, since for $m_{1,n}=s_1\ s_{n}\cdots s_2s_0w_0,$ $\mathcal{L}(m_{1,n})=S-\{s_{n-1},s_2\}$, by Proposition \ref{prop:kl} (2) we see that
$$P_{s_2s_0w_0,m_{1,n}}=P_{m_{1,n-2},m_{1,n}}=1.$$

By direct computations and using (3.4.6)-(3.4.7), we obtain
\begin{equation}
P_{s_0w_0,m_{n-2,n}}=q+1.
\tag{3.4.9}
\end{equation}

Note that we have computed
$$\mu(m_{n-2,n-1},w_{i,n-2}),\ i=2,4,\cdots,n-2$$
in part (6) of the proof of Lemma \ref{lem:B1}, i.e. (3.3.11)-(3.3.12). To complete the calculations in formulas (3.4.1)-(3.4.5), it remains to compute the following Kazhdan-Lusztig polynomials:
\begin{align*}
P_{s_0w_0,m_{n-2,n-1}},\ P_{s_2s_0w_0,w_{3,n-2}},\ P_{s_0w_0,w_{1,n-3}}.
\end{align*}

{\bf (6)} Recall that
$$m_{n-2,n-1}=s_{n-2}\cdots s_1\ s_{n-1}s_{n-2}\cdots s_2s_0w_0.$$
Thus $\mathcal{L}(m_{n-2,n-1})=S-\{s_n,s_0\},$ and by Proposition \ref{prop:kl} (2) we see
\begin{equation}
P_{s_0w_0,m_{n-2,n-1}}=P_{m_{n-2,n-1},m_{n-2,n-1}}=1.
\tag{3.4.10}
\end{equation}

{\bf (7)} Now we compute $P_{s_2s_0w_0,w_{3,n-2}}.$ First we prove some equalities. Recall that
$$w_{3,n-2}=s_3\cdots s_{n-1}s_{n-2}\cdots s_1\ s_{n}\cdots s_2s_0w_0.$$
Thus $\mathcal{L}(w_{3,n-2})=S-\{s_n,s_{2},s_1\}$, and by Proposition \ref{prop:kl} (2) we get
$$P_{s_2s_0w_0,w_{3,n-2}}=P_{m_{n-1},w_{3,n-2}}.$$

For $4\leqslant i\leqslant n-2$,
$$w_{i,n-2}=s_i\cdots s_{n-1}s_{n-2}\cdots s_1\ s_{n}\cdots s_2s_0w_0,$$
and then $w_{i,n-2}=s_{i-1}w_{i-1,n-2}.$ Since
$s_{i-1}m_{n-1}=m_{n-1}s_i$ and $w_{i,n-2}s_i<w_{i,n-2}$. By Proposition \ref{prop:kl} (2) we have
$$P_{s_{i-1}m_{n-1},w_{i,n-2}}=P_{m_{n-1},w_{i,n-2}}.$$
Similarly we have
$$P_{s_{n-2}m_{n-1},m_{n-1,n}}=P_{m_{n-1},m_{n-1,n}}.$$

Keeping in mind these equalities, Using the recursive formula in Proposition \ref{prop:kl} (1) and through  an analysis similar to that for $(*)$, we get the
following three formulas:
\begin{equation}
\begin{aligned}
&P_{m_{n-1},w_{i-1,n-2}}=\\
&\begin{cases}
(1+q)P_{m_{n-1},w_{i,n-2}}-\mu(m_{n-2,n-1},w_{i,n-2})q^{\frac{n-i+2}{2}}P_{m_{n-1},m_{n-2,n-1}}-qP_{m_{n-1},w_{i+1,n-2}},\\
\ \ \ \ \ \ \ \ \ \ \ \ \ \ \ \ \ \ \ \ \ \ \ \ \ \ \ \ \ \ \ \ \ \ \ \ \ \ \ \ \ \ \ \ \ \ \ \ \ \ \ \ \ \ \ \ \ \ \ \ \ \ \ \ \ \ \ \ \ \ \ \
\ \ \ \ \ i=4,6,\cdots,n-4,\\
(1+q)P_{m_{n-1},w_{i,n-2}}-qP_{m_{n-1},w_{i+1,n-2}},
\ \ \ \ \ \ \ \ \ \ \ \ \ \ \ \ \ \ \ \ \ \ \ \ \ \ \ \ \ i=5,7,\cdots,n-3,\\
(1+q)P_{m_{n-1},w_{i,n-2}}-\mu(m_{n-2,n-1},w_{i,n-2})q^{\frac{n-i+2}{2}}P_{m_{n-1},m_{n-2,n-1}}-qP_{m_{n-1},m_{n-1,n}},\\
\ \ \ \ \ \ \ \ \ \ \ \ \ \ \ \ \ \ \ \ \ \ \ \ \ \ \ \ \ \ \ \ \ \ \ \ \ \ \ \ \ \ \ \ \ \ \ \ \ \ \ \ \ \ \ \ \ \ \ \ \ \ \ \ \ \ \ \ \ \ \ \
\ \ \ \ \ i=n-2;
\end{cases}
\nonumber
\end{aligned}
\end{equation}

$$P_{m_{n-1},w_{n-2,n-2}}=(1+q)P_{m_{n-1},m_{n-1,n}}-qP_{m_{n-1},m_{n-2,n}};$$
$$P_{m_{n-1},m_{n-1,n}}=P_{m_{n-1},m_{n-2,n}}+qP_{m_{n-2},m_{n-2,n}}-qP_{m_{n-1},m_{n-2,n-1}}.$$

Since $\mathcal{L}(m_{n-2,n})=S-\{s_{n-1},s_0\}$ by Proposition \ref{prop:kl} (2) we see that
$$P_{m_{n-1},m_{n-2,n}}=P_{m_{n-2,n-1},m_{n-2,n}}=1$$
and (see (3.4.9))
$$P_{m_{n-2},m_{n-2,n}}=P_{s_0w_0,m_{n-2,n}}=q+1.$$
Similarly since $\mathcal{L}(m_{n-2,n-1})=S-\{s_{n},s_0\}$ we see
$$P_{m_{n-1},m_{n-2,n-1}}=P_{m_{n-2,n-1},m_{n-2,n-1}}=1.$$

Combining the previous paragraphs and noting that (3.3.11)-(3.3.12), by a direct computation we obtain
\begin{equation}
P_{s_2s_0w_0,w_{3,n-2}}=q^{n-3}+q+1.
\tag{3.4.11}
\end{equation}

{\bf (8)} By computation and using (3.4.3)-(3.4.5),(3.4.9)-(3.4.11) we get
$$P_{s_0w_0,w_{1,n-2}}=q+1.$$
By (3.4.1)-(3.4.2), we see
$$\deg P_{s_0w_0,z_{3,1}}\leqslant \deg (1+q)^{n-2} P_{s_0w_0,w_{1,n-2}}=n-1,$$
though we don't know $P_{s_0w_0,w_{1,n-3}}$.
Thus
$$\deg P_{s_0w_0,z_{3,1}}\leqslant n-1<2n-4.$$
Since $v_2w_0=z_{3,1}$, the lemma is proved.
\end{proof}

Now we are able to establish the following result, which completes our proof of Theorem \ref{thm:Bmu}.

\begin{Lemma}\label{lem:B3}
Let $d_w$ be as in Theorem \ref{thm:Bmu}. If $n\geqslant 4$ then
$$\mu(w_0,d_ww_0)=1.$$
\end{Lemma}

\begin{proof}
Note that $\mu(w_0,d_ww_0)=\mu(s_0w_0,v_1w_0)$. Then it follows from $(\#)$, Lemma \ref{lem:B1} and Lemma \ref{lem:B2}.
\end{proof}

Next we consider the first extension group between two certain irreducible rational modules of $Sp_{2n}(\bar{\mathbb F}_p)$. Before doing this let's recall some general background
(see \citep{scott2010some}).

We keep notations in Section 2. Assume that $G$ is simply connected and simple and $p$ is a prime number. Denote the root lattice $Q=\mathbb{Z}R.$
Set $s_0'=s_{\alpha_0}t_{p\alpha_0}.$ Let $W'=W_0\ltimes pQ,$ which as an affine transformation group of $E$ is generated by $s_{\alpha}(\alpha\in
R)$ and $s_0'.$ Via these simple reflections, $W_a$ is isomorphic to the group $W'$ by the map: $s_i\mapsto s_i(1\leqslant i\leqslant n),s_0\mapsto s_0'.$ Moreover we
can use these simple reflections to identify $W'$ with the affine Weyl group given by Lusztig in \citep[Section 1.1]{lusztig1980hecke}, i.e., mapping every simple reflection to the orbit of its fixed affine hyperplane. We shall consider that $W'$ acts on the set of weights not in any affine hyperplane $H_{\alpha,k}=\{x\in E|\langle x,\alpha^{\vee}\rangle=pk\}(\alpha\in R,k\in \mathbb{Z})$
through the affine Weyl group  in loc.cit. and denote this action by $*$. Let $\rho$ be the sum of all fundamental weights.
The action $*$ and the original affine action are related:
$w*(-\rho)=w^{-1}(-\rho),w\in W'.$ But for a general element in $E$ similar relation may not hold. In addition, we can easily see that $w\in \Gamma_0$ if and only if
$w*(-\rho)-\rho$ is dominant.

Let $H$ be the algebra group obtained from $G$ by replacing the base field $\mathbb{C}$ with $\bar{\mathbb{F}}_p.$ It is known that when
$p$ is sufficiently large, Lusztig's modular conjecture (see \citep{lusztig1980some}) is true for irreducible modules of $H$
with highest weights in the Jantzen region. The Jantzen region is defined to be the set $\{\nu\in X|0\leqslant \langle
\nu+\rho,\alpha_0^{\vee}\rangle\leqslant p(p-h+2)\},$ where $h$ is the Coxeter number.

Now let $G=Sp_{2n}(\mathbb{C}).$ Let $d_w\in W_a$ be as in Theorem \ref{thm:Bmu} and $\beta=\varepsilon_1+\varepsilon_3.$ Denote the image in $W'$ of $d_w$ under
the above map  by $v'.$ Then $v'=s_{\beta}t_{px_1+px_3}.$ Let $\lambda=t_{2p\rho}w_0*(-\rho)-\rho$ and
$\mu=v't_{2p\rho}w_0*(-\rho)-\rho.$ With the help of the above relation we obtain that
$$\lambda=2p\rho,\ \mu=2p\rho+(p-2n+2)(x_1+x_3)+(2n-2)x_2.$$
Thus $\langle \lambda+\rho,\alpha_0^{\vee}\rangle=p(4n-2)+2n-1$ and $\langle \mu+\rho,\alpha_0^{\vee}\rangle=p(4n+1)+1.$ For
$G=Sp_{2n}(\mathbb{C}),h=2n.$ If $p\geqslant 6n,$ then we have
$$p(p-h+2)\geqslant p(4n+2)>p(4n+1)+1>p(4n-2)+2n-1,$$
i.e. $\lambda$ and $\mu$ are in the Jantzen region. For $H=Sp_{2n}(\bar{\mathbb{F}}_p),$ by Theorem \ref{thm:Bmu} we know
$\mu(t_{2p\rho}w_0,v't_{2p\rho}w_0)=n$ if $n\geqslant 4.$ On the other hand, we know that this   coefficient is identified with the
dimension of the first extension group $\text{Ext}^1_{H}(L(\lambda),L(\mu))$ between the irreducible rational modules $L(\lambda)$ and $L(\mu)$ of $H$ with  highest weights $\lambda$ and
$\mu$  respectively, see for example \citep{andersen1986inversion}. So we obtain the following result.


\begin{Corollary}\label{cor:B}
Let $\lambda$ and $\mu$ be as above and let $L(\lambda)$ and $L(\mu)$ be irreducible rational modules of  $H = Sp_{2n}(\bar{\mathbb{F}}_p)$
with  highest weights $\lambda$ and
$\mu$  respectively.    If $p$ is sufficiently large such that Lusztig's modular conjecture
holds and $p\geq 6n, n\geq 4$, then $\text{Ext}^1_{H}(L(\lambda),L(\mu))=n.$
\end{Corollary}

\section{The case of type $\tilde C_n$ }

In this section we assume that $G=Spin_{2n+1}(\mathbb{C})$ and $W$ is the corresponding extended affine Weyl group. The Coxeter diagram of the affine Weyl group is
\begin{center}
  \begin{tikzpicture}[scale=.6]
    \draw (-1,0) node[anchor=east]  {$\tilde C_n$ };
    \draw[thick] (0 cm,0) circle (.2 cm) node [above] {$0$};
    \draw[thick] (2 cm,0) circle (.2 cm) node [above] {$1$};
    \draw[thick] (4 cm,0) circle (.2 cm) node [above] {$2$};
    \draw[thick] (6 cm,0) circle (.2 cm) node [above] {$n-2$};
    \draw[thick] (8 cm,0) circle (.2 cm) node [above] {$n-1$};
    \draw[thick] (10 cm,0) circle (.2 cm) node [above] {$n$};
    \draw[thick] (.2 cm,.1) -- +(1.6 cm,0);
    \draw[thick] (.2 cm,-.1) -- +(1.6 cm,0);
    \draw[thick] (2.2 cm,0) -- +(1.6 cm,0);
    \draw[dotted,thick] (4.2 cm,0) -- +(1.6 cm,0);
    \draw[thick] (6.2 cm,0) -- +(1.6 cm,0);
    \draw[thick] (8.2 cm, .1 cm) -- +(1.6 cm,0);
    \draw[thick] (8.2 cm, -.1 cm) -- +(1.6 cm,0);
  \end{tikzpicture}
\end{center}
Take a standard orthogonal basis $\varepsilon_i\in E,1\leqslant i\leqslant n,$ such that
$\alpha_1=\varepsilon_1-\varepsilon_2,\alpha_2=\varepsilon_2-\varepsilon_3,\cdots,\alpha_{n-1}=\varepsilon_{n-1}-\varepsilon_{n},
\alpha_n=\varepsilon_n$ are the simple roots and the root system is
$$R=\{\pm(\varepsilon_{i}\pm\varepsilon_{j}),\pm\varepsilon_{k}|1\leqslant i\neq j\leqslant n,1\leqslant k\leqslant n\}.$$
Then $\alpha_0=\varepsilon_1.$ Denote the simple reflections by $s_i=s_{\alpha_i}$ and denote the fundamental dominant weights by
$x_i=x_{\alpha_i},1\leqslant i\leqslant n.$ We have $x_i=\varepsilon_1+\varepsilon_2+\cdots+\varepsilon_i$ for $1\leqslant i<n,$ and
$x_{n}=\frac{1}{2}(\varepsilon_1+\varepsilon_2+\cdots+\varepsilon_{n-1}+\varepsilon_{n}).$

\begin{Proposition}\label{prop:Cdw}
Keep the notations above. If $n\geqslant 4,$ then for
$w=s_{\varepsilon_3}\in W_0,$ we have $x_2<d_w$ and $l(d_w)=l(x_2)+1.$
\end{Proposition}

\begin{proof}
For $w=s_{\varepsilon_3}$, we have $d_w=s_{\varepsilon_3}x_3.$ By computing the lengths, we have $l(d_w)=4n-1$ and
$4n-2=l(x_2)=l(s_2x_2)+1=l(s_1s_2x_2)+2.$ Thus $s_1s_2x_2<s_2x_2<x_2,$ and $l(d_w)=l(x_2)+1.$ It is easy to see $s_2s_1s_0s_1s_2x_2=d_w.$ Hence
$x_2<d_w.$
\end{proof}

Now we can state the main result in this section.

\begin{Theorem}\label{thm:Cmu}
Let $G=Spin_{2n+1}(\mathbb{C}).$ If $n\geqslant 4,$ then for $x=\prod\limits_{i=1}^nx_i^{a_i},a_i\geqslant1,$
$w=s_{\varepsilon_3},$ we have
$$\mu(xw_0,d_wxw_0)=n+1.$$
\end{Theorem}

\begin{proof}
By Lemma \ref{lem:mu} (2), we have
$$\mu(xw_0,d_wxw_0)=\sum\limits_{z_1\in X^+}m_{x^*,x,z_1^*}\delta_{w_0,d_{w} w_0,z_1w_0},$$

Assume that $\delta_{w_0,d_ww_0,z_1w_0}\neq 0$. As the proof of Theorem \ref{thm:Bmu}, we have $z_1\in X^+\cap \mathbb{Z}R,$ and $l(z_1)\leqslant l(d_w)+1=4n.$ By easy calculations, we have
$l(x_1)=2n,l(x_2)=4n-2,l(x_i)=(2n+1-i)i>4n$ for $3\leqslant i<n,l(x_n)=\frac{1}{2}n(n+1),l(x_n^2)=n(n+1)>4n,l(x_1^3)=6n>4n$, and $l(x_1x_2)=6n-2>4n.$ Note
that $x_n\not\in \mathbb{Z}R.$ Therefore $z_1\in \{e,x_1,x_1^2,x_2\}$. So
\begin{equation}
\begin{aligned}
&\mu(xw_0,d_wxw_0)\\
=&m_{x^*,x,e}\delta_{w_0,d_{w} w_0,w_0}+m_{x^*,x,x_1^*}\delta_{w_0,d_{w} w_0,x_1w_0}+m_{x^*,x,(x_1^2)^*}\delta_{w_0,d_{w}
w_0,x_1^2w_0}+m_{x^*,x,x_2^*}\delta_{w_0,d_{w} w_0,x_2w_0}.
\nonumber
\end{aligned}
\end{equation}

Since $l(s_1d_w)=l(d_w)+1=l(x_1^2),s_1x_1^2<x_1^2$ and $s_1d_w>d_w$, so if $d_w<x_1^2$ we must have
$d_w=s_1x_1^2,$
which is impossible by an easy computation. By Lemma \ref{lem:mu} (3) we have $\delta_{w_0,d_ww_0,x_1^2w_0}=\mu(x_1^2w_0,d_ww_0)=0.$ By Proposition \ref{prop:Cdw}, $\delta_{w_0,d_ww_0,x_2w_0}=\mu(x_2w_0,d_ww_0)=1.$

It is easy to see $m_{x^*,x,e}=m_{e,x,x}=1.$ By Lemma \ref{lem:mu} (1) and Proposition \ref{prop:weight}, using Freudenthal's weight multiplicity formula, We obtain
$$m_{x^*,x,x_2^*}=m_{x_2,x,x}=\dim\ V(x_2)_0=n,\ m_{x^*,x,x_1^*}=m_{x_1,x,x}=\dim\ V(x_1)_0=1.$$
Thus
$$\mu(xw_0,d_wxw_0)=\delta_{w_0,d_ww_0,w_0}+\delta_{w_0,d_w w_0,x_1w_0}+n=\mu(w_0,d_{w} w_0)+\mu(x_1w_0,d_{w} w_0)+n.$$

We are now reduced to show that $\mu(w_0,d_{w} w_0)=1$ and $\mu(x_1w_0,d_{w} w_0)=0$, which is Lemma \ref{lem:C3} to be proved next. The theorem is proved.
\end{proof}

 Now we are going to prove $\mu(w_0,d_{w} w_0)=1$ and $\mu(x_1w_0,d_{w} w_0)=0$.
First we need
to work out a reduced expression of $d_w.$ Since $s_0=s_{-\varepsilon_1,1}=s_{\varepsilon_1}t_{\varepsilon_1},$ we have the reduced expression
$$x_1=s_1s_2\cdots s_n\cdots s_2s_1s_0.$$
Thus we get reduced expressions
\begin{equation}
\begin{aligned}
x_2&=s_1x_1s_1x_1\\
&=s_2s_3...s_n...s_3s_2s_1s_0\ s_2s_3...s_n...s_3s_2s_1s_0\\
&=s_2s_1\ s_3...s_n...s_3s_2s_1s_0\ s_3...s_n...s_3s_2s_1s_0,\\
d_w&=s_2s_1s_0s_1s_2x_2\\
&=s_2s_1s_0\ s_3...s_n...s_3s_2s_1s_0\ s_3...s_n...s_3s_2s_1s_0.\nonumber
\end{aligned}
\end{equation}

Set $v=d_w,v_1=s_2d_w.$ Note that $P_{s_2w_0,v_1w_0}=P_{w_0s_2,v_1w_0}=P_{w_0,v_1w_0}$. Using Proposition \ref{prop:kl} (1) we get
\begin{equation}
P_{w_0,vw_0}=(1+q)P_{w_0,v_1w_0}-\sum\limits_{\mbox{\tiny$\begin{array}{c}z\in \Gamma_0\\w_0\leqslant z\prec v_2w_0\\
s_2z<z\end{array}$}}\mu(z,v_1w_0)q^{\frac{1}{2}(l(vw_0)-l(z))}P_{w_0,z}.
\tag{$\clubsuit$}
\end{equation}

We show that the summation in formula $(\clubsuit)$ can be reduced to a single element. That is, there exists only one $z\in \Gamma_0$ such that
$w_0\leqslant z\prec v_1w_0$ and $s_2z<z.$ Assume that $z\in\Gamma_0$ satisfies the conditions.
Since $z\leqslant v_1w_0$, and
$$v_1=s_1s_0\ s_3...s_n...s_3s_2s_1s_0\ s_3...s_n...s_3s_2s_1s_0,$$
by Proposition \ref{prop:leftcell} we see that $z$ must be one of the following elements:
\begin{equation}
\begin{aligned}
&w_0;\\
&s_0w_0;\\
&m_i=s_i\cdots s_1s_0w_0,\ 1\leqslant i\leqslant n;\\
&n_j=s_j\cdots s_{n-1}s_ns_{n-1}\cdots s_1s_0w_0,\ 3\leqslant j\leqslant n-1;\\
&f_i=s_0\ s_i\cdots s_1s_0w_0,\ 1\leqslant i\leqslant n;\\
&g_j=s_0\ s_j\cdots s_{n-1}s_ns_{n-1}\cdots s_1s_0w_0,\ 3\leqslant j\leqslant n-1;\\
&m_{k,i}=s_k\cdots s_1s_0\ \ s_i\cdots s_1s_0w_0,\ 1\leqslant k<i\leqslant n;\\
&n_{k,j}=s_k\cdots s_1s_0\ \ s_j\cdots s_{n-1}s_ns_{n-1}\cdots s_1s_0w_0,\ 1\leqslant k\leqslant n, 3\leqslant j\leqslant n-1;\\
&h_{l,j}=s_l\cdots s_{n-1}s_ns_{n-1}\cdots s_1s_0\ \ s_j\cdots s_{n-1}s_ns_{n-1}\cdots s_1s_0w_0,\ 3\leqslant j\leqslant l\leqslant n-1;\\
&b_{k,i}=s_0\ s_k\cdots s_1s_0\ \ s_i\cdots s_1s_0w_0,\ 1\leqslant k\leqslant n, k< i\leqslant n;\\
&c_{k,j}=s_0\ s_k\cdots s_1s_0\ \ s_j\cdots s_{n-1}s_ns_{n-1}\cdots s_1s_0w_0,\ 1\leqslant k\leqslant n, 3\leqslant j\leqslant n-1;\\
&d_{l,j}=s_0\ s_l\cdots s_{n-1}s_ns_{n-1}\cdots s_1s_0\ \ s_j\cdots s_{n-1}s_ns_{n-1}\cdots s_1s_0w_0,\ 3\leqslant l\leqslant j\leqslant n-1;\\
&u_{k,i}=s_1s_0\ s_k\cdots s_1s_0\ \ s_i\cdots s_1s_0w_0,\ 2\leqslant k<i\leqslant n;\\
&v_{k,j}=s_1s_0\ s_k\cdots s_1s_0\ \ s_j\cdots s_{n-1}s_ns_{n-1}\cdots s_1s_0w_0,\ 2\leqslant k\leqslant n, 3\leqslant j\leqslant n-1;
\\
&w_{l,j}=s_1s_0\ s_l\cdots s_{n-1}s_ns_{n-1}\cdots s_1s_0\ \ s_j\cdots s_{n-1}s_ns_{n-1}\cdots s_1s_0w_0,\ 3\leqslant j\leqslant l\leqslant
n-1;\\
&r_k=s_1\ s_k\cdots s_1s_0\ s_3\cdots s_{n-1}s_ns_{n-1}\cdots s_1s_0w_0,\ 2\leqslant k\leqslant n;\\
&a_{l}=s_1\ s_l\cdots s_{n-1}s_n s_{n-1}\cdots s_1s_0\ s_3\cdots s_{n-1}s_ns_{n-1}\cdots s_1s_0w_0,\ 3\leqslant l\leqslant n-1;\\
&p=s_n\cdots s_1s_0\  s_n\cdots s_1s_0w_0;\\
&t=s_0\ s_n\cdots s_1s_0\  s_n\cdots s_1s_0w_0;\\
&\bar{t}=s_1s_0\ s_n\cdots s_1s_0\  s_n\cdots s_1 s_0w_0;\\
&x=s_2\ s_3\cdots s_{n-1}s_ns_{n-1}\cdots s_1s_0w_0;\\
&\bar{x}=s_1\ s_2\ s_3\cdots s_{n-1}s_ns_{n-1}\cdots s_1s_0w_0;\\
&y=s_2s_0\ s_3\cdots s_{n-1}s_ns_{n-1}\cdots s_1s_0w_0;\\
&\bar{y}=s_1\ s_2s_0\ s_3\cdots s_{n-1}s_ns_{n-1}\cdots s_1s_0w_0.
\nonumber
\end{aligned}
\end{equation}

Since $\mathcal{L}(v_1w_0)=\mathcal{L}(w_{3,3})=S-\{s_0,s_2\}$, combining the condition $s_2z<z$ and Proposition \ref{prop:kl} (4), we see that $z$ is one of
the following elements:
$$w_0,\bar{x},d_{3,3}.$$

Since $l(v_1w_0)-l(w_0)=4n-2$ and $l(v_1w_0)-l(\bar{x})=2n-2$ are even, so $\mu(w_0,v_1w_0)=\mu(\bar{x},v_1w_0)=0.$ Therefore we have
\begin{equation*}
P_{w_0,vw_0}=(1+q)P_{w_0,v_1w_0}-qP_{w_0,d_{3,3}}.
\end{equation*}

Similarly, Note that $d_{3,3}=s_1v_1w_0<v_1w_0,\mathcal{L}(d_{3,3})=S-\{s_1\}$.  Using the recursive formula in Proposition \ref{prop:kl} (1) and through  an analysis similar to that for $(*)$, we get
\begin{equation*}
P_{w_0,v_1w_0}=(1+q)P_{w_0,d_{3,3}}.
\end{equation*}
Thus
\begin{equation*}
P_{w_0,vw_0}=(1+q+q^2)P_{w_0,d_{3,3}}.
\tag{$\heartsuit$}
\end{equation*}

Since $h_{3,3}=s_0d_{3,3}<d_{3,3},$ using the recursive formula in Proposition \ref{prop:kl} (1) and through  an analysis similar to that for $(*)$, noting that $\mathcal{L}(h_{3,3})=S-\{s_0,s_1\},l(h_{3,3})=4n-4+l(w_0),$ we get
\begin{equation}
P_{w_0,d_{3,3}}=P_{w_0,h_{3,3}}+qP_{s_0w_0,h_{3,3}}-\mu(s_0w_0,h_{3,3})q^{2n-2}.
\tag{$\spadesuit$}
\end{equation}



\begin{Lemma}\label{lem:C1}
Let $h_{3,3}$ as above. If $n\geqslant 4$ then
$$
P_{w_0,h_{3,3}}=q^{2n-3}+\text{lower degree terms}.
$$
\end{Lemma}
\begin{proof}
If $n=4,5$ the calculation is easy and we omit it. When $n\geqslant 6$, the case that $n$ is even and the case that $n$
is odd are similar. In view of this, we only consider the case that $n\geqslant 6$ is even in the followings. Recall that we have a reduced expression
$$h_{3,3}w_0=s_3\cdots s_{n-1}s_ns_{n-1}\cdots s_1s_0\ \ s_3\cdots s_{n-1}s_ns_{n-1}\cdots s_1s_0.$$
In the following we shall use an analysis similar to that for $(\clubsuit)$,
to determine the form of $P_{w_0,h_{3,3}}.$
Starting with the above reduced expression, every time we cancel the leftmost simple reflection. We obtain the following (1)-(5).

{\bf (1)} Recall that for $4\leqslant i\leqslant n-1,$
$$h_{i,3}=s_i\cdots s_{n-1}s_ns_{n-1}\cdots s_1s_0\ \ s_3\cdots s_{n-1}s_ns_{n-1}\cdots s_1s_0w_0,$$
and then $h_{i,3}=s_{i-1}h_{i-1,3}.$

Since $\mathcal{L}(h_{i,3})=S-\{s_0,s_1,s_{i-1}\},l(h_{i,3})=(4n-i-1)+l(w_0)$, through  an analysis similar to that for $(\clubsuit)$ we see
that
\begin{equation}
\begin{aligned}
&P_{w_0,h_{i-1,3}}\\
&=
\begin{cases}
(1+q)P_{w_0,h_{i,3}}-\mu(w_0,h_{i,3})q^{2n-2}-\mu(y,h_{i,3})q^{n-2}P_{w_0,y}-q(P_{w_0,h_{4,4}}+P_{w_0,h_{5,3}}),&i=4,\\
(1+q)P_{w_0,h_{i,3}}-\mu(s_0w_0,h_{i,3})q^{2n-\frac{i+1}{2}}-\mu(x,h_{i,3})q^{n-\frac{i-1}{2}}P_{w_0,x}-qP_{w_0,h_{i+1,3}},
&i=5,7,\cdots,n-3,\\
(1+q)P_{w_0,h_{i,3}}-\mu(w_0,h_{i,3})q^{2n-\frac{i}{2}}-\mu(y,h_{i,3})q^{n-\frac{i}{2}}P_{w_0,y}-qP_{w_0,h_{i+1,3}},
&i=6,8,\cdots,n-2,\\
(1+q)P_{w_0,h_{n-1,3}}-\mu(s_0w_0,h_{n-1,3})q^{\frac{3n}{2}}-\mu(x,h_{n-1,3})q^{\frac{n+2}{2}}P_{w_0,x}-qP_{w_0,n_{n,3}},
&i=n-1.
\end{cases}
\end{aligned}
\tag{4.3.1}
\end{equation}
Here
$$x=s_2\ s_3\cdots s_{n-1}s_ns_{n-1}\cdots s_1s_0w_0,$$
$$y=s_2s_0\ s_3\cdots s_{n-1}s_ns_{n-1}\cdots s_1s_0w_0$$
and
$$h_{4,4}=s_4\cdots s_{n-1}s_ns_{n-1}\cdots s_1s_0\ \ s_4\cdots s_{n-1}s_ns_{n-1}\cdots s_1s_0w_0,$$
$$n_{n,3}=s_n\cdots s_1s_0\ \ s_3\cdots s_{n-1}s_ns_{n-1}\cdots s_1s_0w_0.$$

{\bf (2)} Since $n_{n,3}=s_{n-1}h_{n-1,3},\mathcal{L}(n_{n,3})=S-\{s_0,s_1,s_{n-1}\},l(n_{n,3})=(3n-1)+l(w_0),$ through  an analysis similar to that for $(\clubsuit)$ we see that
\begin{equation}
P_{w_0,h_{n-1,3}}=(1+q)P_{w_0,n_{n,3}}-\mu(w_0,n_{n,3})q^{\frac{3n}{2}}-\mu(y,n_{n,3})q^{\frac{n}{2}}P_{w_0,y}-qP_{w_0,n_{n-1,3}}.
\tag{4.3.2}
\end{equation}

Recall that for $2\leqslant i\leqslant n-1,$
$$n_{i,3}=s_i\cdots s_1s_0\ \ s_3\cdots s_{n-1}s_ns_{n-1}\cdots s_1s_0w_0,$$
and then $n_{i,3}=s_{i+1}n_{i+1,3}.$

Since $\mathcal{L}(n_{i,3})=S-\{s_0,s_1,s_{i+1}\},l(n_{i,3})=(2n-1+i)+l(w_0)$, through  an analysis similar to that for $(\clubsuit)$ we see
that
\begin{equation}
\begin{aligned}
&P_{w_0,n_{i+1,3}}\\
&=
\begin{cases}
(1+q)P_{w_0,n_{i,3}}-\mu(s_0w_0,n_{i,3})q^{\frac{2n+i-1}{2}}-\mu(x,n_{i,3})q^{\frac{i+1}{2}}P_{w_0,x}-qP_{w_0,n_{i-1,3}},
&i=n-1,n-3,\cdots,3,\\
(1+q)P_{w_0,n_{i,3}}-\mu(w_0,n_{i,3})q^{\frac{2n+i}{2}}-\mu(y,n_{i,3})q^{\frac{i}{2}}P_{w_0,y}-qP_{w_0,n_{i-1,3}},
&i=n-2,n-4,\cdots,4,\\
(1+q)P_{w_0,n_{i,3}}-\mu(w_0,n_{i,3})q^{n+1}-q(P_{w_0,y}+P_{w_0,n_{i-1,3}}),
&i=2.
\end{cases}
\end{aligned}
\tag{4.3.3}
\end{equation}

Since $n_{1,3}=s_{2}n_{2,3},\mathcal{L}(n_{1,3})=S-\{s_0,s_{2}\},l(n_{1,3})=2n+l(w_0)$, through  an analysis similar to that for $(\clubsuit)$ we see that
\begin{equation}
P_{w_0,n_{2,3}}=(1+q)P_{w_0,n_{1,3}}.
\tag{4.3.4}
\end{equation}

{\bf (3)} Recall that
$$g_3=s_0\ \ s_3\cdots s_{n-1}s_ns_{n-1}\cdots s_1s_0w_0,$$
and then $g_3=s_1n_{1,3}.$

Since $\mathcal{L}(g_3)=S-\{s_1,s_2\},l(g_3)=(2n-1)+l(w_0),$ through  an analysis similar to that for $(\clubsuit)$ we see that
\begin{equation}
P_{w_0,n_{1,3}}=(1+q)P_{w_0,g_3}-qP_{w_0,n_3}.
\tag{4.3.5}
\end{equation}
Here $n_3=s_3\cdots s_{n-1}s_ns_{n-1}\cdots s_1s_0w_0.$

{\bf (4)} Since $n_{3}=s_{0}g_{3},\mathcal{L}(n_3)=S-\{s_0,s_2\},l(n_{3})=(2n-2)+l(w_0)$, we see that
\begin{equation}
P_{w_0,g_3}=P_{w_0,n_{3}}+qP_{s_0w_0,n_{3}}.
\tag{4.3.6}
\end{equation}

Recall that for $4\leqslant i \leqslant n-1,$
$$n_{i}=s_i\cdots s_{n-1}s_ns_{n-1}\cdots s_1s_0w_0,$$
and then $n_{i}=s_{i-1}n_{i-1}.$

Since $\mathcal{L}(n_{i})=S-\{s_0,s_{i-1}\},l(n_{i})=(2n-i+1)+l(w_0)$, through  an analysis similar to that for $(\clubsuit)$ we see that
\begin{equation}
P_{w_0,n_{i-1}}=
\begin{cases}
(1+q)P_{w_0,n_{i}}-\mu(w_0,n_{i})q^{n-\frac{i-2}{2}}-qP_{w_0,n_{i+1}},
&i=4,6,\cdots,n-2,\\
(1+q)P_{w_0,n_{i}}-qP_{w_0,n_{i+1}},
&i=5,7,\cdots,n-3,\\
(1+q)P_{w_0,n_{n-1}}-qP_{w_0,m_n},
&i=n-1.
\end{cases}
\tag{4.3.7}
\end{equation}
Here $m_n=s_ns_{n-1}\cdots s_1s_0w_0.$

{\bf (5)} Since $m_n=s_{n-1}n_{n-1},\mathcal{L}(m_n)=S-\{s_0,s_{n-1}\},l(m_n)=(n+1)+l(w_0)$, we see that
\begin{equation}
P_{w_0,n_{n-1}}=(1+q)P_{w_0,m_n}-\mu(w_0,m_n)q^{\frac{n}{2}+1}-qP_{w_0,m_{n-1}}.
\tag{4.3.8}
\end{equation}

{\bf (6)} through  an analysis similar to that for $(\clubsuit)$ we can easily get that
$$P_{w_0,m_{i}}=1,$$
for $1\leqslant i\leqslant n.$
Thus by direct computation with (4.3.7)-(4.3.8), we see that
\begin{equation}
P_{w_0,n_{i}}=1,3\leqslant i\leqslant n-1.
\tag{4.3.9}
\end{equation}
Similarly we can easily get that
\begin{equation}
P_{s_0w_0,n_{i}}=1,3\leqslant i\leqslant n-1.
\tag{4.3.10}
\end{equation}
Thus by (4.3.6), we obtain
\begin{equation}
P_{w_0,g_3}=q+1.
\tag{4.3.11}
\end{equation}
Hence for completing the calculation in (4.3.1)-(4.3.5), it remains to know the followings:
\begin{align*}
&\mu(y,h_{i,3}),i=4,6,\cdots,n-2;\\
&\mu(y,n_{i,3}),i=n,n-2,\cdots,4;\\
&\mu(x,h_{i,3}),i=5,7,\cdots,n-1;\\
&\mu(x,n_{i,3}),i=n-1,n-3,\cdots,3;\\
&\mu(s_0w_0,h_{i,3}),i=5,7,\cdots,n-1;\\
&\mu(s_0w_0,n_{i,3}),i=n-1,n-3,\cdots,3;\\
&P_{w_0,y},P_{w_0,x};\\
&P_{w_0,h_{4,4}}.
\end{align*}

{\bf (7)} Recall that for $4\leqslant i\leqslant n-1$,
$$h_{i,3}=s_i\cdots s_{n-1}s_ns_{n-1}\cdots s_1s_0\ \ s_3\cdots s_{n-1}s_ns_{n-1}\cdots s_1s_0w_0;$$
for $3\leqslant i\leqslant n,$
$$n_{i,3}=s_i\cdots s_1s_0\ \ s_3\cdots s_{n-1}s_ns_{n-1}\cdots s_1s_0w_0;$$
and
$$y=s_2s_0\ s_3\cdots s_{n-1}s_ns_{n-1}\cdots s_1s_0w_0.$$

For $4\leqslant i\leqslant n$, $s_iy=ys_i$, thus
$$P_{s_iy,h_{j,3}}=P_{y,h_{j,3}},5\leqslant j\leqslant n-1,$$
and
$$P_{s_iy,n_{k,3}}=P_{y,n_{k,3}},3\leqslant k\leqslant n.$$

Therefore through  an analysis similar to that for $(\clubsuit)$ and referring to (4.3.1)-(4.3.3), we see (note that $x<y$)
\begin{align*}
&P_{y,h_{i-1,3}}=
\begin{cases}
(1+q)P_{y,h_{i,3}}-qP_{y,h_{i+1,3}},
&i=5,7,\cdots,n-3,\\
(1+q)P_{y,h_{i,3}}-\mu(y,h_{i,3})q^{n-\frac{i}{2}}-qP_{y,h_{i+1,3}},
&i=6,8,\cdots,n-2,\\
(1+q)P_{y,h_{n-1,3}}-qP_{y,n_{n,3}},
&i=n-1.
\end{cases}\\
&P_{y,h_{n-1,3}}=(1+q)P_{y,n_{n,3}}-\mu(y,n_{n,3})q^{\frac{n}{2}}-qP_{y,n_{n-1,3}}.\\
&P_{y,n_{i+1,3}}=
\begin{cases}
(1+q)P_{y,n_{i,3}}-qP_{y,n_{i-1,3}},
&i=n-1,n-3,\cdots,3,\\
(1+q)P_{y,n_{i,3}}-\mu(y,n_{i,3})q^{\frac{i}{2}}-qP_{y,n_{i-1,3}},
&i=n-2,n-4,\cdots,4.
\end{cases}
\end{align*}

Since $l(n_{3,3})-2=l(n_{2,3})-1=l(y),$ we have
$$P_{y,n_{3,3}}=P_{y,n_{2,3}}=1.$$
By direct computation, we see
\begin{equation}
\begin{aligned}
&P_{y,h_{j,3}}=1,\mu(y,h_{j,3})=0,4\leqslant j\leqslant n-1;\\
&P_{y,n_{k,3}}=1,\mu(y,n_{k,3})=0,4\leqslant k\leqslant n.
\end{aligned}
\tag{4.3.12}
\end{equation}

Recall that
$$x=s_2\ s_3\cdots s_{n-1}s_ns_{n-1}\cdots s_1s_0w_0.$$
Using a similar argument, we can obtain that
\begin{equation}
\begin{aligned}
&P_{x,h_{j,3}}=1,\mu(x,h_{j,3})=0,5\leqslant j\leqslant n-1;\\
&P_{x,n_{k,3}}=1,\mu(x,n_{k,3})=0,3\leqslant k\leqslant n-1.
\end{aligned}
\tag{4.3.13}
\end{equation}

{\bf (8)} Note that
$$g_3=s_2y=s_0\ \ s_3\cdots s_{n-1}s_ns_{n-1}\cdots s_1s_0w_0.$$
Since $\mathcal{L}(g_3)=S-\{s_1,s_2\},l(g_3)=(2n-1)+l(w_0),$ through  an analysis similar to that for $(\clubsuit)$ we see that
$$P_{w_0,y}=(1+q)P_{w_0,g_3}-qP_{w_0,g_4}.$$
Here $g_4=s_0\ \ s_4\cdots s_{n-1}s_ns_{n-1}\cdots s_1s_0w_0.$

Since $n_4=s_0g_4,\mathcal{L}(n_4)=S-\{s_3,s_0\},l(n_4)=(2n-3)+l(w_0),$ through  an analysis similar to that for $(\clubsuit)$ we see that
$$P_{w_0,g_4}=P_{w_0,n_4}+qP_{s_0w_0,n_4}.$$

By (4.3.9)-(4.3.11), we can get
\begin{equation}
P_{w_0,y}=q+1.
\tag{4.3.14}
\end{equation}

{\bf (9)} Recall that
$$h_{4,4}=s_4\cdots s_{n-1}s_ns_{n-1}\cdots s_1s_0\ \ s_4\cdots s_{n-1}s_ns_{n-1}\cdots s_1s_0w_0.$$

Since $h_{4,4}=s_2h_{4,3},\mathcal{L}(h_{4,4})=S-\{s_2,s_0\},l(h_{4,4})=(4n-6)+l(w_0),$ through  an analysis similar to that for $(\clubsuit)$ we see that
$$P_{w_0,h_{4,3}}=(1+q)h_{4,4}.$$

From Lemma \ref{lem:C2} below  we will see that (see (4.4.11))
\begin{align*}
&\mu(s_0w_0,h_{i,3})=0,i=5,7,\cdots,n-1,\\
&\mu(s_0w_0,n_{i,3})=0,i=n-1,n-3,\cdots,3,
\end{align*}
which together with (4.3.1)-(4.3.14) imply
\begin{equation}
\begin{aligned}
&P_{w_0,h_{4,3}}=q^{2n-3}+2q^{2n-4}+2q^{2n-5}+q^{2n-6}+\cdots+q+1;\\
&P_{w_0,h_{5,3}}=q^{2n-4}+2q^{2n-5}+2q^{2n-6}+q^{2n-7}+\cdots+q+1.
\end{aligned}
\tag{4.3.15}
\end{equation}
Thus
\begin{equation}
P_{w_0,h_{4,4}}=q^{2n-4}+q^{2n-5}+q^{2n-6}+q^{2n-8}+\cdots+q^2+1.
\tag{4.3.16}
\end{equation}

{\bf (10)} Finally by (4.3.1)
\begin{align*}
P_{w_0,h_{3,3}}=(1+q)P_{w_0,h_{4,3}}-\mu(w_0,h_{4,3})q^{2n-2}-\mu(y,h_{4,3})q^{n-2}P_{w_0,y}-q(P_{w_0,h_{4,4}}+P_{w_0,h_{5,3}}).
\end{align*}
By (4.3.15)-(4.3.16) and (4.3.12) we get
$$P_{w_0,h_{3,3}}=q^{2n-3}+\text{lower degree terms}.$$
\end{proof}

\begin{Lemma}\label{lem:C2}
Let $h_{3,3}$ be as in Lemma \ref{lem:C1}. If $n\geqslant 4$ then
$$P_{s_0w_0,h_{3,3}}=q^{2n-3}+1.$$
\end{Lemma}
\begin{proof}
As Lemma \ref{lem:C1}, we shall only consider the case that $n\geqslant 6$ is even
in the following. Using an analysis similar to that for $(\clubsuit)$, we obtain the following (1)-(4) (refer to the corresponding ones in Lemma \ref{lem:C1} and note
specially that $g_3=s_1n_{1,3},n_3=s_0g_3$).

{\bf (1)}
\begin{equation}
\begin{aligned}
&P_{s_0w_0,h_{i-1,3}}\\
&=
\begin{cases}
(1+q)P_{s_0w_0,h_{4,3}}-\mu(y,h_{4,3})q^{n-2}P_{s_0w_0,y}-q(P_{s_0w_0,h_{4,4}}+P_{s_0w_0,h_{5,3}}),
&i=4,\\
(1+q)P_{s_0w_0,h_{i,3}}-\mu(s_0w_0,h_{i,3})q^{2n-\frac{i+1}{2}}-\mu(x,h_{i,3})q^{n-\frac{i-1}{2}}P_{s_0w_0,x}-qP_{s_0w_0,h_{i+1,3}},
&i=5,7,\cdots,n-3,\\
(1+q)P_{s_0w_0,h_{i,3}}-\mu(y,h_{i,3})q^{n-\frac{i}{2}}P_{s_0w_0,y}-qP_{s_0w_0,h_{i+1,3}},&i=6,8,\cdots,n-2,\\
(1+q)P_{s_0w_0,h_{n-1,3}}-\mu(s_0w_0,h_{n-1,3})q^{\frac{3n}{2}}-\mu(x,h_{n-1,3})q^{\frac{n+2}{2}}P_{s_0w_0,x}-qP_{s_0w_0,n_{n,3}},
&i=n-1.
\end{cases}
\end{aligned}
\tag{4.4.1}
\end{equation}

{\bf (2)}
\begin{equation}
P_{s_0w_0,h_{n-1,3}}=(1+q)P_{s_0w_0,n_{n,3}}-\mu(y,n_{n,3})q^{\frac{n}{2}}P_{s_0w_0,y}-qP_{s_0w_0,n_{n-1,3}}.
\tag{4.4.2}
\end{equation}

\begin{equation}
\begin{aligned}
&P_{s_0w_0,n_{i+1,3}}\\
&=
\begin{cases}
(1+q)P_{s_0w_0,n_{i,3}}-\mu(s_0w_0,n_{i,3})q^{\frac{2n+i-1}{2}}-\mu(x,n_{i,3})q^{\frac{i+1}{2}}P_{s_0w_0,x}-qP_{s_0w_0,n_{i-1,3}},
&i=n-1,n-3,\cdots,3,\\
(1+q)P_{s_0w_0,n_{i,3}}-\mu(y,n_{i,3})q^{\frac{i}{2}}P_{s_0w_0,y}-qP_{s_0w_0,n_{i-1,3}},&i=n-2,n-4,\cdots,4,\\
(1+q)P_{s_0w_0,n_{i,3}}-q(P_{s_0w_0,y}+P_{s_0w_0,n_{i-1,3}}),&i=2.
\end{cases}
\end{aligned}
\tag{4.4.3}
\end{equation}

\begin{equation}
P_{s_0w_0,n_{2,3}}=(1+q)P_{s_0w_0,n_{1,3}}.
\tag{4.4.4}
\end{equation}

{\bf (3)}
\begin{equation}
P_{s_0w_0,n_{1,3}}=P_{s_0w_0,g_3}+qP_{s_1s_0w_0,g_3}-qP_{s_0w_0,n_3}.
\tag{4.4.5}
\end{equation}

{\bf (4)}
\begin{equation}
P_{s_0w_0,g_3}=P_{w_0,n_{3}}+qP_{s_0w_0,n_{3}}.
\tag{4.4.6}
\end{equation}

{\bf (5)} Recall that
$$n_3=s_3\cdots s_{n-1}s_ns_{n-1}\cdots s_1s_0w_0.$$
Thus by Proposition 1.7 we have $s_0s_1s_0w_0\not\leqslant n_3$.

Since
$n_3=s_0g_3,\mathcal{L}(n_3)=S-\{s_0,s_2\},l(n_3)=(2n-2)+l(w_0),$ through an analysis similar to that for $(\clubsuit)$ we see that
\begin{equation}
P_{s_1s_0w_0,g_3}=P_{s_1s_0w_0,n_{3}}.
\tag{4.4.7}
\end{equation}

Easily we get $P_{s_1s_0w_0,n_{3}}=1$ (refer to the process of getting (4.3.9)). Thus by (4.3.9)-(4.3.11) and (4.4.5)-(4.4.7) we obtain
\begin{equation}
P_{s_0w_0,n_{1,3}}=q+1.
\tag{4.4.8}
\end{equation}

Hence for completing the calculation in (4.4.1)-(4.4.4), it remains to compute the followings:
$$P_{s_0w_0,y},
P_{s_0w_0,h_{4,4}}.
$$

{\bf (6)} Since $s_0y<y,$ noting (4.3.14), we have
\begin{equation}
P_{s_0w_0,y}=P_{w_0,y}=q+1.
\tag{4.4.9}
\end{equation}

{\bf (7)} Using (4.4.1)-(4.4.4) and (4.4.8) to do calculation, we get
\begin{equation}
\begin{aligned}
&P_{s_0w_0,h_{5,3}}=q^{2n-5}+1,\\
&P_{s_0w_0,h_{4,3}}=(1+q)(q^{2n-5}+1),
\end{aligned}
\tag{4.4.10}
\end{equation}
and
\begin{equation}
\begin{aligned}
&\mu(s_0w_0,h_{i,3})=0,i=5,7,\cdots,n-1,\\
&\mu(s_0w_0,n_{i,3})=0,i=n-1,n-3,\cdots,3.
\end{aligned}
\tag{4.4.11}
\end{equation}

Since $h_{4,4}=s_2h_{4,3},$ through an analysis similar to that for $(\clubsuit)$ we see that
$$P_{s_0w_0,h_{4,3}}=(1+q)P_{s_0w_0,h_{4,4}}.$$
Thus
\begin{equation}
P_{s_0w_0,h_{4,4}}=q^{2n-5}+1,
\tag{4.4.12}
\end{equation}

{\bf (8)} Finally by (4.4.1)
$$P_{s_0w_0,h_{3,3}}=(1+q)P_{s_0w_0,h_{4,3}}-\mu(y,h_{4,3})q^{n-2}P_{s_0w_0,y}-q(P_{s_0w_0,h_{4,4}}+P_{s_0w_0,h_{5,3}}).$$
By (4.3.12),(4.4.9),(4.4.10),(4.4.12), we get
$$P_{s_0w_0,h_{3,3}}=q^{2n-3}+1.$$
\end{proof}

\begin{Lemma}\label{lem:C3}
Let $d_w$ be as in Theorem \ref{thm:Cmu}. If $n\geqslant 4$ then
$$\mu(w_0,d_ww_0)=1,\ \mu(x_1w_0,d_ww_0)=0.$$
\end{Lemma}

\begin{proof}
The first equality follows from $(\heartsuit),(\spadesuit)$, Lemma \ref{lem:C1} and Lemma \ref{lem:C2}. Now we see the second. As above we only
consider the case that $n\geqslant 6$ is even.

Using an analysis similar to that for $(\clubsuit)$, we have
$$P_{x_1w_0,vw_0}=(1+q)P_{x_1w_0,v_1w_0}.$$
Since $d_{3,3}=s_1v_1w_0<v_1w_0$ and
$x_1w_0\not\leqslant d_{3,3}$ by Proposition \ref{prop:kl}(3) we obtain
$$P_{x_1w_0,v_1w_0}=P_{x,d_{3,3}}.$$

Referring to $(\spadesuit)$ and (4.3.1), we get (note that $y=s_0x$)
\begin{align*}
&P_{x,d_{3,3}}=P_{x,h_{3,3}}+qP_{y,h_{3,3}},\\
&P_{x,h_{3,3}}=(1+q)P_{x,h_{4,3}}-\mu(y,h_{4,3})q^{n-2}-qP_{x,h_{5,3}},\\
&P_{x,h_{4,3}}=(1+q)P_{x,h_{5,3}}-\mu(x,h_{5,3})q^{n-2}-qP_{x,h_{6,3}},\\
&P_{y,h_{3,3}}=(1+q)P_{y,h_{4,3}}-\mu(y,h_{4,3})q^{n-2}-qP_{y,h_{5,3}},
\end{align*}
which together with (4.3.12),(4.3.13) imply that $P_{x,d_{3,3}}=q+1.$ Thus $P_{x_1w_0,d_ww_0}=q^2+2q+1.$
Hence
$$\mu(x_1w_0,d_ww_0)=0.$$
\end{proof}


 Next we consider the first extension group between two certain irreducible rational modules of $Spin_{2n+1}(\bar{ \mathbb{F}}_p)$, where $p$ is a prime number.

Now let $G=Spin_{2n+1}(\mathbb{C})$. Let $d_w$ be as in Theorem \ref{thm:Cmu}. Denote the root lattice $Q=\mathbb{Z}R.$ Let $W^\prime = W_0\ltimes pQ$, which is isomorphic to the affine Weyl group $W_a$.  Denote the image in $W'$ of
$d_w$ under the isomorphism by $v'$. Then $v^\prime = s_{\varepsilon_3}t_{px_3}.$  Let $\lambda=t_{2p\rho}w_0*(-\rho)-\rho$ and $\mu=v't_{2p\rho}w_0*(-\rho)-\rho.$ By a simple computation, we get
\[
\lambda = 2p\rho,\ \mu = 2p\rho + (p-2n+5)x_3 + (2n-5)x_2.
\]
Thus $\langle \lambda+\rho, \alpha_0^\vee\rangle=p(4n-2)+(2n-1)$ and $\langle \mu+\rho, \alpha_0^\vee\rangle=4np+(2n-1)$. For $G=Spin_{2n+1}(\mathbb{C})$, the Coxeter number is $h=2n$. If
$p\geq 6n-1$, then we have
\[
p(p-h+2)\geqslant p(4n+1) > 4np+(2n-1) > p(4n-2)+(2n-1).
\]
Therefore $\lambda$ and $\mu$ are in the Jantzen region. For $H=Spin_{2n+1}(\bar{\mathbb{F}}_p)$, by Theorem \ref{thm:Cmu} we know $\mu(t_{p\rho}w_0,v^\prime t_{p\rho}w_0)=n+1$ if
$n\geqslant 4$. By the relation of the leading coefficient and the dimension of the first extension group, we get the following result.

\begin{Corollary}\label{cor:C}
Let $\lambda$ and $\mu$ be as above and let $L(\lambda)$ and $L(\mu)$ be irreducible rational modules of  $H = Spin_{2n+1}(\bar{\mathbb{F}}_p)$  with  highest weights $\lambda$ and
$\mu$  respectively.    If $p$ is sufficiently large such that Lusztig's modular conjecture
holds and $p\geq 6n-1,n\geq 4$, then $\text{Ext}^1_{H}(L(\lambda),L(\mu))=n+1.$
\end{Corollary}

\section{The case of type $\tilde D_n$ }

In this section we assume that $G=Spin_{2n}(\mathbb{C})$ and $W$ is the corresponding extended affine Weyl group.  The Coxeter diagram of the affine Weyl group is
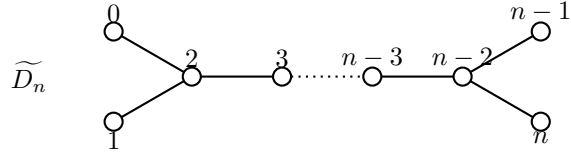
\begin{center}
  \begin{tikzpicture}[scale=.6]
    \draw (-1,0) node[anchor=east]  {$\tilde D_n$ };
    \draw[thick] (2 cm,0) circle (.2 cm) node [above] {$2$};
    \draw[xshift=2 cm,thick] (150:2) circle (.2 cm) node [above] {$0$};
    \draw[xshift=2 cm,thick] (210:2) circle (.2 cm) node [below] {$1$};
    \draw[thick] (4 cm,0) circle (.2 cm) node [above] {$3$};
    \draw[thick] (6 cm,0) circle (.2 cm) node [above] {$n-3$};
    \draw[thick] (8 cm,0) circle (.2 cm) node [above] {$n-2$};
    \draw[xshift=8 cm,thick] (30:2) circle (.2 cm) node [above] {$n-1$};
    \draw[xshift=8 cm,thick] (-30:2) circle (.2 cm) node [below] {$n$};
    \draw[xshift=2 cm,thick] (150:0.2) -- (150:1.8);
    \draw[xshift=2 cm,thick] (210:0.2) -- (210:1.8);
    \draw[thick] (2.2,0) --+ (1.6,0);
    \draw[dotted,thick] (4.2,0) --+ (1.6,0);
    \draw[thick] (6.2,0) --+ (1.6,0);
    \draw[xshift=8 cm,thick] (30:0.2) -- (30:1.8);
    \draw[xshift=8 cm,thick] (-30:0.2) -- (-30:1.8);
  \end{tikzpicture}
\end{center}
Take a standard orthogonal basis $\varepsilon_i\in E,1\leqslant i\leqslant n,$ such that
$\alpha_1=\varepsilon_1-\varepsilon_2,\alpha_2=\varepsilon_2-\varepsilon_3,\cdots,\alpha_{n-1}=\varepsilon_{n-1}-\varepsilon_{n},
\alpha_n=\varepsilon_{n-1}+\varepsilon_n$ are the simple roots and the root system
$$R=\{\pm(\varepsilon_{i}\pm\varepsilon_{j})|1\leqslant i\neq j\leqslant n\}.$$
Then $\alpha_0=\varepsilon_1+\varepsilon_2.$ Denote the simple reflections by $s_i=s_{\alpha_i}$ and denote the fundamental dominant weights by
$x_i=x_{\alpha_i},1\leqslant i\leqslant n.$ We have that $x_i=\varepsilon_1+\varepsilon_2+\cdots+\varepsilon_i,$ for $1\leqslant i<n-1,$ and
$x_{n-1}=\frac{1}{2}(\varepsilon_1+\varepsilon_2+\cdots+\varepsilon_{n-1}-\varepsilon_{n}),$
$x_{n}=\frac{1}{2}(\varepsilon_1+\varepsilon_2+\cdots+\varepsilon_{n-1}+\varepsilon_{n}).$

\begin{Proposition}\label{prop:Ddw}
Keep the notations above. If $n\geqslant 5$ then for
$w=s_{\varepsilon_1+\varepsilon_3}\in W_0,$ we have $x_2<d_w$ and $l(d_w)=l(x_2)+1.$
\end{Proposition}

\begin{proof}
For $d_w=s_{\varepsilon_1+\varepsilon_3}$, we have $d_w=s_{\varepsilon_1+\varepsilon_3}x_1x_3=s_2s_0s_2x_2$. By computing the lengths, we have
$l(d_w)=l(x_2)+1=l(s_2x_2)+2=4n-5.$ Thus $s_2x_2<x_2<d_w$ and $l(d_w)=l(x_2)+1.$
\end{proof}

\begin{Theorem}\label{thm:Dmu}
Let $G=Spin_{2n}(\mathbb{C}).$ If $n\geqslant 5,$ then for $x=\prod\limits_{i=1}^nx_i^{a_i},a_i\geqslant1,$
$w=s_{\varepsilon_1+\varepsilon_3},$ we have
$$\mu(xw_0,d_wxw_0)=n+1.$$
\end{Theorem}

\begin{proof}
By Lemma \ref{lem:mu} (2) we have
$$\mu(xw_0,d_wxw_0)=\sum\limits_{z_1\in X^+}m_{x^*,x,z_1^*}\delta_{w_0,d_{w} w_0,z_1w_0}.$$

Assume that $\delta_{w_0,d_ww_0,z_1w_0}\neq 0$. As the proof of Theorem \ref{thm:Bmu}, we have $z_1\in X^+\cap \mathbb{Z}R,$ and $l(z_1)\leqslant l(d_w)+1=4n-4.$ By easy calculations, we have
$l(x_1)=2n-2,l(x_2)=4n-6,l(x_i)=(2n-1-i)i>4n-4$ for $3\leqslant
i<n-1,l(x_{n-1})=l(x_n)=\frac{1}{2}n(n-1),l(x_{n-1}^2)=l(x_{n-1}x_n)=l(x_n^2)=n(n-1)>4n-4,l(x_1^3)=6n-6>4n-4$, and $l(x_1x_2)=6n-8>4n-4.$ Note that
$x_1,x_{n-1},x_n\not\in \mathbb{Z}R.$ Therefore  $z_1\in \{e,x_1^2,x_2\}$. So
\begin{equation}
\begin{aligned}
&\mu(xw_0,d_wxw_0)\\
=&m_{x^*,x,e}\delta_{w_0,d_{w} w_0,w_0}+m_{x^*,x,(x_1^2)^*}\delta_{w_0,d_{w} w_0,x_1^2w_0}+m_{x^*,x,x_2^*}\delta_{w_0,d_{w} w_0,x_2w_0}.
\nonumber
\end{aligned}
\end{equation}

Since $l(s_1d_w)=l(d_w)+1=l(x_1^2),s_1x_1<x_1$, so if $d_w<x_1^2$ we must have $s_1d_w=x_1^2,$ which is impossible by an easy computation. By Lemma \ref{lem:mu} (3) we have
$\delta_{w_0,d_ww_0,x_1^2w_0} = \mu(x_1^2w_0,d_ww_0)=0.$ By Proposition \ref{prop:Ddw}, $\delta_{w_0,d_ww_0,x_2w_0}=\mu(x_2w_0,d_ww_0)=1.$

By Lemma \ref{lem:mu} (1) and Proposition \ref{prop:weight}, using Freudenthal's multiplicity formula, We obtain $m_{x^*,x,x_2^*}=m_{x_2,x,x}=\dim\
V(x_2)_0=n.$ Thus
$$\mu(xw_0,d_wxw_0)=\mu(w_0,d_{w} w_0)+n.$$

We are now reduced to show that $\mu(w_0,d_{w} w_0)=1$, which is Lemma \ref{lem:D3} to be proved next. The theorem is proved.
\end{proof}

 Now we are going to prove $\mu(w_0,d_ww_0)=1.$ First we need to get a reduced expression of $d_w.$ Since
$s_0=s_{\alpha_0}t_{\alpha_0}=s_{\varepsilon_1+\varepsilon_2}x_2$, so
$$x_2=s_{\varepsilon_1+\varepsilon_2}s_0=s_2\cdots s_{n-1}\ s_1\cdots s_{n-2}s_ns_{n-2}\cdots s_1\ s_{n-1}\cdots s_2s_0.$$
Thus we get a reduced expression
\begin{align*}
d_w&=s_2s_0s_2x_2\\
&=s_2s_0\ s_3\cdots s_{n-1}\ s_1\cdots s_{n-2}s_ns_{n-2}\cdots s_1\ s_{n-1}\cdots s_2s_0.
\end{align*}

Set $v=d_w,v_1=s_2d_w$ and $v_2=s_0s_2d_w.$
Since $w_0,vw_0\in \mathscr{D}_L(s_0,s_2)=\{w\in W_a|\#\mathcal{L}(w)\cap \{s_0,s_2\}=1\}$, using the star operator, we have
$\mu(w_0,vw_0)=\mu(s_0w_0,v_1w_0).$ Using Proposition  \ref{prop:kl} (1) we get
\begin{equation}
P_{s_0w_0,v_1w_0}=P_{w_0,v_2w_0}+qP_{s_0w_0,v_2w_0}-\sum\limits_{\mbox{\tiny$\begin{array}{c}z\in \Gamma_0\\s_0w_0\leqslant z\prec v_2w_0\\
s_0z<z\end{array}$}}\mu(z,v_2w_0)q^{\frac{1}{2}(l(v_1w_0)-l(z))}P_{s_0w_0,z}.
\tag{$\diamondsuit$}
\end{equation}

We show that the summation in formula $(\diamondsuit)$ is empty. That is, there exists no $z\in \Gamma_0$ such that $s_0w_0\leqslant z\prec v_2w_0$ and
$s_0z<z$. Assume that $z=uw_0$, for some $u\in W$ with $l(z)=l(u)+l(w_0),$ satisfies the conditions. Since $s_0w_0\leqslant z\leqslant
v_2w_0,$ and
$$v_2=s_3\cdots s_{n-1}\ s_1s_2\cdots s_{n-2}s_{n}s_{n-2}\cdots s_2s_1\ s_{n-1} \cdots s_3s_2\ s_0,$$
by Proposition \ref{prop:leftcell} we see that $z$ must be one of the following elements:
\begin{equation}
\begin{aligned}
&s_0w_0;\\
&m_i=s_i\cdots s_3s_2s_0w_0,\ 2\leqslant i\leqslant n-1;\\
&m_{j,i}=s_j\cdots s_2s_1\ \ s_i\cdots s_2s_0w_0,\ 1\leqslant j<i\leqslant n-1;\\
&f_j=s_ns_j\cdots s_1\ s_{n-2}\cdots s_2s_0w_0,\ 1\leqslant j\leqslant n-3;\\
&g_j=s_ns_j\cdots s_1\ s_{n-1}\cdots s_2s_0w_0,\ 1\leqslant j\leqslant n-2;\\
&h_k=s_k\cdots s_{n-2}s_n\ \ s_{n-1}\cdots s_2s_0w_0,\ 1\leqslant k\leqslant n-2;\\
&w_{k,j}=s_k\cdots s_{n-2}s_n s_j\cdots s_2s_1\ \ s_{n-1}\cdots s_2s_0w_0,\ 1\leqslant k\leqslant n-2, 1\leqslant j\leqslant n-2;\\
&b_{k}=s_{n-1}\ \ s_k\cdots s_{n-2} s_ns_{n-3}\cdots s_1\ \ s_{n-1}\cdots s_2s_0w_0,\ 1\leqslant k\leqslant n-2;\\
&y_{t,k}=s_t\cdots s_{n-1}\ \ s_k\cdots s_{n-2} s_ns_{n-2}\cdots s_1\ \ s_{n-1}\cdots s_2s_0w_0,\ 3\leqslant t\leqslant n-1,1\leqslant k<t;\\
&p=s_n\ s_{n-2}\cdots s_2s_0w_0;\\
&\bar{p}=s_n\ s_{n-1}\cdots s_2s_0w_0;\\
&r=s_{n-2}s_ns_{n-3}\cdots s_1\ s_{n-2}\cdots s_2s_0w_0;\\
&\bar{r}=s_{n-1}\ s_{n-2}s_ns_{n-3}\cdots s_1\ s_{n-2}\cdots s_2s_0w_0.
\nonumber
\end{aligned}
\end{equation}

Since $\mathcal{L}(v_2w_0)=\mathcal{L}(y_{3,1})=S-\{s_2,s_0\},$ combining the condition $s_0z<z$ and Proposition \ref{prop:kl} (4), we see that $z$ only can be $s_0w_0$ or $w_{3,1}.$ Since
$l(v_2w_0)-l(s_0w_0)=4n-6$ and $l(v_2w_0)-l(w_{3,1})=2n-4$, so $\mu(s_0w_0,v_2w_0)=\mu(w_{3,1},v_2w_0)=0.$ Therefore we have
\begin{equation}
P_{s_0w_0,v_1w_0}=P_{w_0,v_2w_0}+qP_{s_0w_0,v_2w_0}.
\tag{$\blacklozenge$}
\end{equation}

\begin{Lemma}\label{lem:D1}
Let $v_2$ be as above. If $n\geqslant 5$ then
$$P_{w_0,v_2w_0}=q^{2n-4}+\text{lower degree terms}.$$
\end{Lemma}
\begin{proof}
If $n=5,6$ the calculation is easy and we omit it. When $n\geqslant 7$, the case that $n$ is even and the case that $n$
is odd are similar. In view of this, we only consider the case that $n\geqslant 7$ is odd in the followings.
Recall that we have reduced expression
$$v_2=s_3\cdots s_{n-1}\ s_1\cdots s_{n-2}s_ns_{n-2}\cdots s_1\ s_{n-1}\cdots s_2s_0.$$
In the following we shall use repeatedly an analysis similar to that for $(\diamondsuit)$
to determine the form of $P_{w_0,v_2w_0}.$
Starting with the above reduced expression, every time we cancel the leftmost simple reflection. We obtain the following (1)-(4). Note that $v_2w_0=y_{3,1}.$

{\bf (1)} Recall that for $4\leqslant i\leqslant n-1$,
$$y_{i,1}=s_i\cdots s_{n-1}\ \ s_1\cdots s_{n-2} s_ns_{n-2}\cdots s_1\ \ s_{n-1}\cdots s_2s_0w_0,$$
and then $y_{i,1}=s_{i-1}y_{i-1,1}.$

Since $\mathcal{L}(y_{i,1})=S-\{s_{i-1},s_0\},l(y_{i,1})=(4n-4-i)+l(w_0)$, through  an analysis similar to that for $(\diamondsuit)$ we see
\begin{equation}
P_{w_0,y_{i-1,1}}=
\begin{cases}
(1+q)P_{w_0,y_{i,1}}-qP_{w_0,y_{i+1,1}},
&i=4,6,\cdots,n-3,\\
(1+q)P_{w_0,y_{i,1}}-\mu(w_0,y_{i,1})q^{2n-2-\frac{i-1}{2}}-qP_{w_0,y_{i+1,1}},&i=5,7,\cdots,n-2,\\
(1+q)P_{w_0,y_{n-1,1}}-q(P_{w_0,b_1}+P_{w_0,w_{1,n-2}}),&i=n-1.
\end{cases}
\tag{5.3.1}
\end{equation}
Here
$$b_1=s_{n-1}\ \ s_1\cdots s_{n-2} s_ns_{n-3}\cdots s_1\ \ s_{n-1}\cdots s_2s_0w_0$$
and
$$w_{1,n-2}=s_1\cdots s_{n-2}s_n s_{n-2}\cdots s_2s_1\ \ s_{n-1}\cdots s_2s_0w_0.$$

{\bf (2)} Since $w_{1,n-2}=s_{n-1}y_{n-1,1},\mathcal{L}(w_{1,n-2})=S-\{s_{n-1},s_0\},l(w_{1,n-2})=(3n-4)+l(w_0),$ through  an analysis similar to that for $(\diamondsuit)$ we see that
\begin{equation}
P_{w_0,y_{n-1,1}}=(1+q)P_{w_0,w_{1,n-2}}-\mu(w_0,w_{1,n-2})q^{\frac{3n-3}{2}}.
\tag{5.3.2}
\end{equation}

Recall that for $2\leqslant i\leqslant n-2$,
$$w_{i,n-2}=s_i\cdots s_{n-2}s_n s_{n-2}\cdots s_2s_1\ \ s_{n-1}\cdots s_2s_0w_0,$$
and then $w_{i,n-2}=s_{i-1}w_{i-1,n-2}.$

Since $\mathcal{L}(w_{i,n-2})=S-\{s_{i-1},s_{n-1},s_0\},l(w_{i,n-2})=(3n-3-i)+l(w_0)$, through  an analysis similar to that for $(\diamondsuit)$ we see that
\begin{equation}
P_{w_0,w_{i-1,n-2}}=
\begin{cases}
(1+q)P_{w_0,w_{i,n-2}}-\mu(r,w_{i,n-2})q^{\frac{n-i+1}{2}}-qP_{w_0,w_{i+1,n-2}},&i=2,4,\cdots,n-3,\\
(1+q)P_{w_0,w_{i,n-2}}-\mu(w_0,w_{i,n-2})q^{\frac{3n-2-i}{2}}-qP_{w_0,w_{i+1,n-2}},&i=3,5,\cdots,n-4,\\
(1+q)P_{w_0,w_{n-2,n-2}}-\mu(w_0,w_{n-2,n-2})q^{n}-qP_{w_0,g_{n-2}},&i=n-2.
\end{cases}
\tag{5.3.3}
\end{equation}
Here
$$r=s_{n-2}s_ns_{n-3}\cdots s_1\ s_{n-2}\cdots s_2s_0w_0$$
and
$$g_{n-2}=s_n s_{n-2}\cdots s_2s_1\ \ s_{n-1}\cdots s_2s_0w_0.$$

{\bf (3)} Since $g_{n-2}=s_{n-2}w_{n-2,n-2},\mathcal{L}(g_{n-2})=S-\{s_{n-2},s_0\},l(g_{n-2})=(2n-2)+l(w_0),$ through  an analysis similar to that for $(\diamondsuit)$ we see that
\begin{equation}
P_{w_0,w_{n-2,n-2}}=(1+q)P_{w_0,g_{n-2}}-qP_{w_0,m_{n-2,n-1}}.
\tag{5.3.4}
\end{equation}
Here $m_{n-2,n-1}=s_{n-2}\cdots s_2s_1\ \ s_{n-1}\cdots s_2s_0w_0.$

{\bf (4)} Since $m_{n-2,n-1}=s_ng_{n-2},\mathcal{L}(m_{n-2,n-1})=S-\{s_n,s_0\},l(m_{n-2,n-1})=(2n-3)+l(w_0),$ through  an analysis similar to that for $(\diamondsuit)$ we see that
\begin{equation}
P_{w_0,g_{n-2}}=(1+q)P_{w_0,m_{n-2,n-1}}-\mu(w_0,m_{n-2,n-1})q^{n-1}.
\tag{5.3.5}
\end{equation}

{\bf (5)} In fact, we can use the result about the case of $\widetilde{B_n}$ to obtain $P_{w_0,m_{n-2,n-1}}.$ Note that we have the sub-diagram
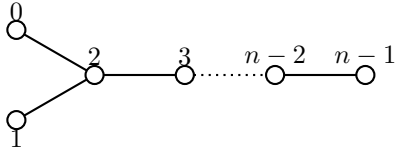
\begin{center}
  \begin{tikzpicture}[scale=.6]
    \draw[thick] (2 cm,0) circle (.2 cm) node [above] {$2$};
    \draw[xshift=2 cm,thick] (150:2) circle (.2 cm) node [above] {$0$};
    \draw[xshift=2 cm,thick] (210:2) circle (.2 cm) node [below] {$1$};
    \draw[thick] (4 cm,0) circle (.2 cm) node [above] {$3$};
    \draw[thick] (6 cm,0) circle (.2 cm) node [above] {$n-2$};
    \draw[thick] (8 cm,0) circle (.2 cm) node [above] {$n-1$};
    \draw[xshift=2 cm,thick] (150:0.2) -- (150:1.8);
    \draw[xshift=2 cm,thick] (210:0.2) -- (210:1.8);
    \draw[thick] (2.2,0) --+ (1.6,0);
    \draw[dotted,thick] (4.2,0) --+ (1.6,0);
    \draw[thick] (6.2,0) --+ (1.6,0);
  \end{tikzpicture}
\end{center}
which is same as the corresponding sub-diagram of $\widetilde{B_n}.$ If we note the parity of $n$ then by (3.3.15) we obtain that
\begin{equation*}
P_{w_0,m_{n-2,n-1}}=q^{n-3}+q^{n-5}+\cdots+q^2+1.
\tag{5.3.6}
\end{equation*}

Hence for completing the calculation in (1)-(4), it remains to compute the followings:
\begin{align*}
&\mu(r,w_{i,n-2}),i=2,4,\cdots,n-3;\\
&P_{w_0,r};\\
&P_{w_0,b_1}.
\end{align*}

{\bf (6)} Now we compute $P_{r,w_{2,n-2}}.$ Recall that
$$r=s_{n-2}s_ns_{n-3}\cdots s_1\ s_{n-2}\cdots s_2s_0w_0.$$
For $2\leqslant
i\leqslant n-3,$ noting that $s_jw_0=w_0s_{j}$
($1\leqslant j\leqslant n-2$) and $s_nw_0=w_0s_{n-1}$, we have $s_ir=rs_{i+2}.$
Thus
$$P_{r,w_{i+1,n-2}}=P_{s_ir,w_{i+1,n-2}} \text{, for } 2\leqslant i\leqslant n-3.$$
Therefore through  an analysis similar to that for $(\diamondsuit)$ and referring to (5.3.3), we see that
\begin{equation*}
P_{r,w_{i-1,n-2}}=
\begin{cases}
(1+q)P_{r,w_{i,n-2}}-\mu(r,w_{i,n-2})q^{\frac{n-i+1}{2}}-qP_{r,w_{i+1,n-2}},&i=4,6\cdots,n-3,\\
(1+q)P_{r,w_{i,n-2}}-qP_{r,w_{i+1,n-2}},&i=3,5,\cdots,n-4,\\
(1+q)P_{r,w_{n-2,n-2}},&i=n-2.
\end{cases}
\end{equation*}

Since $l(w_{n-2,n-2})=l(r)+2$, we have
$$P_{r,w_{n-2,n-2}}=1.$$
By computation with the above equalities, we obtain
$$P_{r,w_{i,n-2}}=q+1,2\leqslant i\leqslant n-3,$$
thus
\begin{equation}
\mu(r,w_{i,n-2})=
\begin{cases}
0,&i=2,4,\cdots,n-5,\\
1,&i=n-3.
\end{cases}
\tag{5.3.7}
\end{equation}

{\bf (7)} First we introduce the following natural isomorphism of groups:
\begin{align*}
\sigma:W_a&\rightarrow W_a\\
s_i&\mapsto s_i, 0\leqslant i\leqslant n-2\\
s_{n-1}&\mapsto s_n\\
s_n&\mapsto s_{n-1}.
\end{align*}
It is obvious that this isomorphism preserves the Bruhat order and the length function. Thus the Kazhdan-Lusztig polynomials keep stable under
the action of $\sigma$. Note that $w_0$ is the unique longest element in $W_0$, so it is invariant under this isomorphism.

Recall that
\begin{align*}
&r=s_{n-2}s_ns_{n-3}\cdots s_1\ s_{n-2}\cdots s_2s_0w_0,\\
&m_{n-2,n-1}=s_{n-2}s_{n-3}\cdots s_1\ s_{n-1}\cdots s_2s_0w_0\\
&{\phantom{m_{n-2,n-1}}}=s_{n-2}s_{n-1}s_{n-3}\cdots s_1\ s_{n-2}\cdots s_2s_0w_0.
\end{align*}
Thus $r=\sigma(m_{n-2,n-1})$ and $P_{w_0,r}=P_{w_0,m_{n-2,n-1}}$.
Therefore by (5.3.6) we have
\begin{equation}
P_{w_0,r}=q^{n-3}+q^{n-5}+\cdots+q^2+1.
\tag{5.3.8}
\end{equation}

{\bf (8)} By (5.3.2)-(5.3.7), we see that
\begin{equation}
\begin{aligned}
&P_{w_0,w_{1,n-2}}=q^{n-1}+2q^{n-2}+q^{n-3}+q^{n-4}+\cdots+q+1,\\
&P_{w_0,y_{n-1,1}}=(1+q)(q^{n-1}+2q^{n-2}+q^{n-3}+q^{n-4}+\cdots+q+1).
\end{aligned}
\tag{5.3.9}
\end{equation}

On the other hand, noting that $b_1=s_ny_{n-1,1}$ and through  an analysis similar to that for $(\diamondsuit)$ we see that
$$P_{w_0,y_{n-1,1}}=(1+q)P_{w_0,b_1}-\mu(w_0,b_1)q^{\frac{3n-3}{2}}.$$

Since $1+q|P_{w_0,y_{n-1,1}},$ so $\mu(w_0,b_1)=0.$ Hence
\begin{equation}
P_{w_0,b_1}=q^{n-1}+2q^{n-2}+q^{n-3}+q^{n-4}+\cdots+q+1.
\tag{5.3.10}
\end{equation}

{\bf (9)} By (5.3.1),(5.3.9) and (5.3.10), we see
$$P_{w_0,y_{n-2,1}}=q^{n+1}+\text{lower degree terms}.$$
It is of strictly larger degree than $P_{w_0,y_{n-1,1}}.$

Proceeding with (5.3.1), we find
$$\mu(w_0,y_{i+1,1})=0,i=4,6,\cdots,n-3,$$
thus
$$P_{w_0,v_{2}w_0}=P_{w_0,y_{3,1}}=q^{2n-4}+\text{lower degree terms}.$$
\end{proof}

\begin{Lemma}\label{lem:D2}
Let $v_2$ be as in Lemma \ref{lem:D1}. If $n\geqslant 5$ then
$$\deg P_{s_0w_0,v_2w_0}<2n-5.$$
\end{Lemma}
\begin{proof}
As Lemma \ref{lem:D1}, we shall only consider the case that $n\geqslant 7$ is odd in
the followings. Using  an analysis similar to that for $(\diamondsuit)$ repeatly, we obtain the following (1)-(4) (refer to the corresponding ones in Lemma \ref{lem:D1} and note
specially that $w_{3,n-2}=s_2w_{2,n-2}$).

{\bf (1)}
\begin{equation}
P_{s_0w_0,y_{i-1,1}}=
\begin{cases}
(1+q)P_{s_0w_0,y_{i,1}}-qP_{s_0w_0,y_{i+1,1}},&4\leqslant i\leqslant n-2,\\
(1+q)P_{s_0w_0,y_{n-1,1}}-q(P_{s_0w_0,b_1}+P_{s_0w_0,w_{1,n-2}},&i=n-1.
\end{cases}
\tag{5.4.1}
\end{equation}

{\bf (2)}
\begin{equation}
P_{s_0w_0,y_{n-1,1}}=(1+q)P_{s_0w_0,w_{1,n-2}}.
\tag{5.4.2}
\end{equation}

\begin{equation}
P_{s_0w_0,w_{i-1,n-2}}=
\begin{cases}
(1+q)P_{s_0w_0,w_{i,n-2}}-\mu(r,w_{i,n-2})q^{\frac{n-i+1}{2}}P_{s_0w_0,r}-qP_{s_0w_0,w_{i+1,n-2}},&i=2,4,\cdots,n-3,\\
P_{s_0w_0,w_{i,n-2}}+qP_{s_2s_0w_0,w_{i,n-2}}-qP_{s_0w_0,w_{i+1,n-2}},&i=3,\\
(1+q)P_{s_0w_0,w_{i,n-2}}-qP_{s_0w_0,w_{i+1,n-2}},&i=5,\cdots,n-4,\\
(1+q)P_{s_0w_0,w_{n-2,n-2}}-qP_{s_0w_0,g_{n-2}},&i=n-2.
\end{cases}
\tag{5.4.3}
\end{equation}

{\bf (3)}
\begin{equation}
P_{s_0w_0,w_{n-2,n-2}}=(1+q)P_{s_0w_0,g_{n-2}}-qP_{s_0w_0,m_{n-2,n-1}}.
\tag{5.4.4}
\end{equation}

{\bf (4)}
\begin{equation}
P_{s_0w_0,g_{n-2}}=(1+q)P_{s_0w_0,m_{n-2,n-1}}.
\tag{5.4.5}
\end{equation}

{\bf (5)} Recall that $m_{n-2,n-1}=s_{n-2}\cdots s_1\ s_{n-1}\cdots s_2s_0w_0,\mathcal{L}(m_{n-2,n-1})=S-\{s_n,s_0\}$. Thus
\begin{equation}
P_{s_0w_0,m_{n-2,n-1}}=P_{m_{n-2,n-1},m_{n-2,n-1}}=1.
\tag{5.4.6}
\end{equation}

Similarly, since $r=s_{n-2}s_ns_{n-3}\cdots s_1\ s_{n-2}\cdots s_2s_0w_0$ and $\mathcal{L}(r)=S-\{s_{n-1},s_0\}$, we have
\begin{equation}
P_{s_0w_0,r}=P_{r,r}=1.
\tag{5.4.7}
\end{equation}

{\bf (6)} Recall that
$$f_1=s_ns_1\ s_{n-2}\cdots s_2s_0w_0,$$
$$w_{3,n-2}=s_3\cdots s_{n-2}s_n s_{n-2}\cdots s_2s_1\ \ s_{n-1}\cdots s_2s_0w_0$$
and
$\mathcal{L}(w_{3,n-2})=S-\{s_2,s_{n-1},s_0\}.$
Thus $P_{s_2s_0w_0,w_{3,n-2}}=P_{f_1,w_{3,n-2}}.$

Since $s_if_1=f_1s_{i+1}$ for $3\leqslant i\leqslant n-2,$ Using  an analysis similar to that for $(\diamondsuit)$ and referring to (5.4.3)-(5.4.4) we see that
\begin{equation*}
P_{f_1,w_{i-1,n-2}}=
\begin{cases}
(1+q)P_{f_1,w_{i,n-2}}-\mu(r,w_{i,n-2})q^{\frac{n-i+1}{2}}P_{f_1,r}-qP_{f_1w_{i+1,n-2}},&i=4,6,\cdots,n-3,\\
(1+q)P_{f_1,w_{i,n-2}}-qP_{f_1,w_{i+1,n-2}},&i=5,7,\cdots,n-4,\\
(1+q)P_{f_1,w_{n-2,n-2}}-qP_{f_1,g_{n-2}},&i=n-2;
\end{cases}
\end{equation*}
and
\begin{equation*}
P_{f_1,w_{n-2,n-2}}=(1+q)P_{f_1,g_{n-2}}.
\end{equation*}

As (5.4.7) we have
$$P_{f_1,r}=P_{r,r}=1.$$

Recall that
$$f_{n-3}=s_n\ s_{n-3}\cdots s_2s_1\ s_{n-2}\cdots s_2s_0w_0,$$
$$g_{n-2}=s_n\ s_{n-2}\cdots s_2s_1\ s_{n-1}\cdots s_2s_0w_0$$
and
$\mathcal{L}(g_{n-2})=S-\{s_{n-2},s_0\}.$
Thus $P_{f_1,g_{n-2}}=P_{f_{n-3},g_{n-2}}=1,$ from $l(g_{n-2})-l(f_{n-3})=2.$

Recall that (5.3.7) says that
\begin{equation*}
\mu(r,w_{i,n-2})=
\begin{cases}
0,&i=2,4,\cdots,n-5,\\
1,&i=n-3,
\end{cases}
\end{equation*}
which together with the above computation implies that
$$P_{f_1,w_{i,n-2}}=q^{n-1-i}+q+1,3\leqslant i\leqslant n-3.$$
Thus
\begin{equation}
P_{s_2s_0w_0,w_{3,n-2}}=q^{n-4}+q+1.
\tag{5.4.8}
\end{equation}

{\bf (7)} By (5.4.3)-(5.4.7) and (5.3.7), we get
\begin{equation}
\begin{aligned}
&P_{s_0w_0,w_{4,n-2}}=q^{n-4}+q^{n-5}+q+1,\\
&P_{s_0w_0,w_{3,n-2}}=q^{n-3}+q^{n-4}+q+1.
\end{aligned}
\tag{5.4.9}
\end{equation}

Note that (5.4.3) tells us
\begin{align*}
&P_{s_0w_0,w_{1,n-2}}=P_{s_0w_0,w_{2,n-2}}-qP_{s_0w_0,w_{3,n-2}},\\
&P_{s_0w_0,w_{2,n-2}}=P_{s_0w_0,w_{3,n-2}}+qP_{s_2s_0w_0,w_{3,n-2}}-qP_{s_0w_0,w_{4,n-2}},
\end{align*}
which together with (5.4.8)-(5.4.9) imply that
$$\deg P_{s_0w_0,w_{1,n-2}}\leqslant n-3.$$
Then from (5.4.1)-(5.4.3) we can easily obtain
$$\deg P_{s_0w_0,y_{3,1}}\leqslant \deg P_{s_0w_0,w_{1,n-2}}+n-3<2n-5.$$
\end{proof}

\begin{Lemma}\label{lem:D3}
Let $d_w$ be as in Theorem \ref{thm:Dmu}. If $n\geqslant 5$ then
$$\mu(w_0,d_ww_0)=1.$$
\end{Lemma}

\begin{proof}
Note that $\mu(w_0,d_ww_0)=\mu(s_0w_0,v_1w_0)$ and $l(v_1)=4n-6.$ It follows from $(\blacklozenge)$, Lemma \ref{lem:D1} and Lemma \ref{lem:D2}.
\end{proof}

Next we consider the first extension group between two certain irreducible rational modules of $Spin_{2n}(\bar{ \mathbb{F}}_p)$, where $p$ is a prime number.

Now let $G=Spin_{2n}(\mathbb{C})$. Let $d_w$ be as in Theorem \ref{thm:Dmu}. Denote the root lattice $Q=\mathbb{Z}R.$ Let $W^\prime = W_0\ltimes pQ$, which is isomorphic to the affine Weyl group $W_a$.  Denote the image in $W'$ of
$d_w$ under the isomorphism by $v'$. Then $v^\prime = s_{\varepsilon_1+\varepsilon_3}t_{px_1+px_3}.$  Let $\lambda=t_{2p\rho}w_0*(-\rho)-\rho$ and $\mu=v't_{2p\rho}w_0*(-\rho)-\rho.$ By a simple computation, we get
\[
\lambda = 2p\rho,\ \mu = 2p\rho + (p-2n+4)(x_1+x_3) + (2n-4)x_2.
\]
Thus $\langle \lambda+\rho, \alpha_0^\vee\rangle=p(4n-6)+(2n-3)$ and $\langle \mu+\rho, \alpha_0^\vee\rangle=p(4n-3)+1$. For $G=Spin_{2n}(\mathbb{C})$, the Coxeter number is $h=2n-2$. If
$p\geqslant 6n-6$, then we have
\[
p(p-h+2)\geqslant p(4n-2)> p(4n-3)+1 > p(4n-6)+(2n-3).
\]
Therefore $\lambda$ and $\mu$ are in the Jantzen region. For $H=Spin_{2n}(\bar{\mathbb{F}}_p)$, by Theorem \ref{thm:Dmu} we know $\mu(t_{p\rho}w_0,v^\prime t_{p\rho}w_0)=n+1$ if
$n\geqslant 5$. By the relation of the leading coefficient and the dimension of the first extension group, we get the following result.

\begin{Corollary}\label{cor:D}
Let $\lambda$ and $\mu$ be as above and let $L(\lambda)$ and $L(\mu)$ be irreducible rational modules of  $H = Spin_{2n}(\bar{\mathbb{F}}_p)$
with  highest weights $\lambda$ and
$\mu$  respectively.    If $p$ is sufficiently large such that Lusztig's modular conjecture
holds and $p\geq 6n-6, n\geq 5$, then $\text{Ext}^1_{H}(L(\lambda),L(\mu))=n+1.$
\end{Corollary}




\end{document}